\documentclass[12pt]{amsart}
\usepackage{amssymb}
\usepackage{mathrsfs, amssymb, array}
\usepackage{verbatim}
\usepackage[matrix,arrow,curve]{xy}
\usepackage{longtable}
\usepackage{hhline}
\usepackage{array}
\usepackage{enumitem}
\setenumerate[0]{leftmargin=23pt}

\newcommand{\ZZ}{\mathbb{Z}}
\newcommand{\PP}{\mathbb{P}}
\newcommand{\QQ}{\mathbb{Q}}
\newcommand{\FF}{\mathbb{F}}
\newcommand{\RR}{\mathbb{R}}
\newcommand{\OOO}{\mathscr{O}}
\newcommand{\LLL}{\mathscr{L}}

\newcommand{\Sym}{\mathfrak{S}}
\newcommand{\Alt}{\mathfrak{A}}
\newcommand{\DD}{\mathfrak{D}}

\newcommand{\rk}{\operatorname{rk}}
\newcommand{\Sing}{\operatorname{Sing}}

\newcommand{\Bir}{\operatorname{Bir}}

\newcommand{\Pic}{\operatorname{Pic}}

\newcommand{\Aut}{\operatorname{Aut}}

\newcommand{\Cr}{\operatorname{Cr}}

\newcommand{\tr}{\operatorname{tr}}
\newcommand{\sgn}{\operatorname{sgn}}

\newcommand{\ord}{\operatorname{ord}}
\newcommand{\diag}{\operatorname{diag}}

\newcommand{\id}{\operatorname{id}}

\newcommand{\mQ}{\operatorname{Q}}
\newcommand{\mW}{\operatorname{W}}

\newcommand{\Eu}{\operatorname{Eu}}

\newcommand{\xref}[1]{{\rm\ref{#1}}}

\newcommand{\comp}{\mathbin{\scriptstyle{\circ}}}

\newtheorem{theorem}{Theorem}
\newtheorem{stheorem}[equation]{}

\numberwithin{theorem}{section}
\numberwithin{equation}{theorem}

\newtheorem{mtheorem}[theorem]{}

\theoremstyle{definition}
\newtheorem{say}[theorem]{}

\newtheorem{ppar}[equation]{}
\theoremstyle{definition}

\title{On stable conjugacy of finite subgroups of the plane Cremona group, II}

\author{Yuri Prokhorov}
\address{
Steklov Mathematical Institute, 
8 Gubkina street, Moscow 119991
\newline\indent
Department 
of Algebra, Moscow State
University, Moscow 117234
\newline\indent
Laboratory of Algebraic Geom., SU-HSE, 
7 Vavilova str., Moscow 117312
}

\email{prokhoro@gmail.com}

\thanks{
Partially supported by RScF grant no. 14-21-00053.
}

\begin{document}
\maketitle
 \begin{abstract}
We prove that, except for a few cases,
stable linearizability of finite subgroups of 
the plane Cremona group implies linearizability.
 \end{abstract}

\section{Introduction}
This is a follow-up paper to \cite{Bogomolov-Prokhorov}.
Let $\Bbbk$ be an algebraically closed field of characteristic $0$.
Recall that the \textit{Cremona group} $\Cr_n(\Bbbk)$  is 
the group birational automorphisms $\Bir(\PP^n)$ of the projective space $\PP^n$ over $\Bbbk$.
Subgroups $G\subset \Cr_n(\Bbbk)$ and $G'\subset \Cr_{m}(\Bbbk)$ 
are said to be \textit{stably conjugate} if, for some $N\ge n,\, m$, they
are conjugate in $\Cr_{N}(\Bbbk)$,
where 
the embeddings $\Cr_{n}(\Bbbk),\, \Cr_{m}(\Bbbk)\subset \Cr_{N}(\Bbbk)$ 
are induced by birational isomorphisms 
$\PP^N \dashrightarrow \PP^{n}\times\PP^{N-n}\dashrightarrow \PP^{m}\times\PP^{N-m}$.



Any embedding of a finite subgroup $G\subset \Cr_n(\Bbbk)$ is induced by
a biregular action on a  rational  variety $X$.
A subgroup $G\subset \Cr_n(\Bbbk)$ is said to be \textit{linearizable} 
if one can take $X=\PP^n$. A subgroup $G\subset \Cr_n(\Bbbk)$ is said to be
\textit{stably linearizable} if it is stably conjugate to a linear
action of $G$ on a vector space $\Bbbk^m$.

The following question is a natural extension of the famous Zariski cancellation
problem \cite{BCTSSD85} to the geometric situation.

\begin{mtheorem}{\bf Question.}\label{Question}
Let $G\subset \Cr_2(\Bbbk)$ be a stably linearizable finite subgroup.
It it true that $G$ is linearizable?
\end{mtheorem}

In the present paper we give a partial answer to this question.
In fact, we find a (very restrictive) list of all subgroups $G\subset \Cr_2(\Bbbk)$
which potentially can give counterexamples to \ref{Question}.

It is easy to show (see \cite{Bogomolov-Prokhorov}) that the group $H^1(G,\Pic(X))$ 
is stable birational invariant.
In particular, if $G\subset \Cr_n(\Bbbk)$ is stably linearizable, 
then $H^1(G_1,\Pic(X))=0$ for any subgroup $G_1\subset G$
(then we say that $G\subset \Cr_n(\Bbbk)$ is \textit{$H^1$-trivial}).
Any finite subgroup $G\subset \Cr_2(\Bbbk)$ 
is induced by an action on either a del Pezzo surface or a conic 
bundle \cite{Iskovskikh-1979s-e}.
In the first case our main result is the following theorem
which is based on a computation of $H^1(G,\Pic(X))$ in \cite{Bogomolov-Prokhorov}
(see Theorem \ref{theorem-main-p}).

\begin{mtheorem}{\bf Theorem.}\label{Theorem-main-del-Pezzo}
Let $X$ be a del Pezzo surface
and let $G\subset \Aut(X)$ be a finite subgroup such that the pair
$(X,G)$ is minimal. Then the following are equivalent:
\begin{enumerate}
\item \label{Theorem-main-del-Pezzo-1}
$H^1(G_1,\Pic(X))=0$ for any subgroup $G_1\subset G$,
\item\label{Theorem-main-del-Pezzo-2}
any element of $G$ does not fix a curve of positive genus,
\item\label{Theorem-main-del-Pezzo-3}
either
\begin{enumerate}[leftmargin=17pt]
\item
$K_X^2\ge 5$, or 
\item\hspace{-5pt}\footnote{This case is missing in \cite[Th. 6.9]{Dolgachev-Iskovskikh}.
This is because the arguments on page 489 (case 3) are incorrect.
However, $X$ has an equivariant rational curve fibration (see Remark \ref{Remark-d=4-case-realized}). So, 
the description of the group appears in \cite[Th. 5.7]{Dolgachev-Iskovskikh}.
Note that groups $(\ZZ/2\ZZ)^2{}_{\bullet}\Sym_3$ and $(\ZZ/2\ZZ)^3{}_{\bullet}\Sym_3$
are also missing in \cite[Th. 6.9]{Dolgachev-Iskovskikh}.
}\quad
$X$ is a quartic del Pezzo surface given by 
\begin{equation}\label{equation-d=4}
x^2_1 + \zeta_3 x^2_2 +\zeta_3^2 x^2_3 + x^2_4= x^2_1 + \zeta_3^2 x^2_2 + \zeta_3 x^2_3 + x^2_5 = 0,
\end{equation}
where $\zeta_3=\exp(2\pi i/3)$ and $G\simeq (\ZZ/3\ZZ)\rtimes (\ZZ/4\ZZ)$
is generated by the following two transformations:
\begin{equation}\label{equation-d=4-action}
\begin{array}{lll}
\gamma: (x_1,x_2,x_3,x_4, x_5) &\longmapsto& (x_2,x_3,x_1,\zeta_3x_4,\zeta_3^2x_5), 
\\[5pt]
\beta': (x_1,x_2,x_3,x_4, x_5) &\longmapsto& (x_1, x_3,x_2,-x_5,x_4).
\end{array}
\end{equation}
\end{enumerate}\end{enumerate}
\end{mtheorem}

The conic bundle case is considered in \S \ref{section-Conic-bundles}.
Main results are Theorems 
\ref{theorem-conic-bundle-1} and \ref{theorem-conic-bundle-2}.

Note that there are only a few subgroups $G\subset \Cr_2(\Bbbk)$ that 
are not linearizable and
satisfy the equivalent conditions \ref{Theorem-main-del-Pezzo-1}-\ref{Theorem-main-del-Pezzo-3}
above (see \cite[\S 8]{Dolgachev-Iskovskikh}).

The plan of the proof of Theorem \ref{Theorem-main-del-Pezzo}
is the following.
The most difficult part of the proof is 
the implication \ref{Theorem-main-del-Pezzo-2}$\Rightarrow$\ref{Theorem-main-del-Pezzo-3}.
It is proved in \S\ref{section-d=4}-\S\ref{section-d=1}.
The implication \ref{Theorem-main-del-Pezzo-1}$\Rightarrow$\ref{Theorem-main-del-Pezzo-2}
is exactly the statement of Corollary \ref{Corollary-H1}
and \ref{Theorem-main-del-Pezzo-3}$\Rightarrow$\ref{Theorem-main-del-Pezzo-1}
is a consequence of Proposition \ref{lemma-H1-trivial-d>5} and Corollary \ref{Corollary-lemma-H1-trivial-d=4}.

We tried to make the paper self-contained as much as possible, so 
in the proofs we do not use detailed lists from the classification of
finite subgroups of $\Cr_2(\Bbbk)$ \cite{Dolgachev-Iskovskikh}.
Instead of this we tried to use just \emph{general} facts and principles of this classification. 

\par\medskip\noindent
{\bf Acknowledgments.}
I like to thank Igor Dolgachev for useful conversations through
e-mails. I am also grateful to the referee for constructive criticism  and for pointing me
out a gap in the earlier version of Lemma \ref{Corollary-cyclic}.

\section{Preliminaries}
\begin{say}{\bf Notation.}\label{notation-1}
\begin{itemize}
\item
$\Sym_n$ is the symmetric group.
\item
$\sgn: \Sym_n\to \{\pm 1\}$ is the sign map.
\item
$\mathfrak A_n$ is the alternating group.
\item
$\DD_n$ is a dihedral group of order $2n$,
$n\ge 2$ (in particular, $\DD_2\simeq (\ZZ/2\ZZ)^{2}$).
We will use the following presentation
\begin{equation}\label{equation-dihedral-group}
\DD_n=\langle r,\, s \mid r^n=s^2=1,\ srs=r^{-1}\rangle.
\end{equation}
\item
$\sigma: \DD_n\to \{\pm 1\}$ is the homomorphism defined by 
$\sigma(r)=1$, $\sigma(s)=-1$.
\item
$\tilde \DD_{n}$ is the binary dihedral group (see e.g. \cite{Springer1977}). We identify $\tilde \DD_{n}$ with
the subgroup of $SL_2(\Bbbk)$ generated by the matrices
\begin{equation}\label{equation-binary-dihedral-group}
\tilde r = 
\begin{pmatrix}
\zeta_{2n}& 0
\\
0& \zeta_{2n}^{-1}
\end{pmatrix}
\qquad
\tilde s =
\begin{pmatrix}
0_{\phantom{n}}& i_{\phantom{n}}
\\
i_{\phantom{n}}& 0_{\phantom{n}}
\end{pmatrix}
\end{equation}
Note that $\tilde \DD_{n}$ is
a non-trivial central extension of $\DD_{n}$ 
by $\ZZ/2\ZZ$.
\item
$\zeta_n$ is a primitive $n$-th root of unity.
\item
$\Phi_n(t)$ is the $n$-th cyclotomic polynomial.
\item
$\Eu(X)$ is the topological Euler number of $X$.
\item
$\diag(a_1,\dots,a_n)$ is the diagonal matrix.
\item $X^G$ is the fixed point locus of an action of $G$ on $X$.
\end{itemize}

\end{say}
\begin{say}{\bf $G$-varieties.}
Throughout this paper $G$ denotes a finite group.
We use the standard language of $G$-varieties (see e.g. \cite{Dolgachev-Iskovskikh}).
In particular, we systematically use the following fact:
for any projective non-singular $G$-surface $X$ there exists 
a birational $G$-equivariant morphism $X\to X_{\min}$ such that the $G$-surface $X_{\min}$ 
is \emph{$G$-minimal}, that is, any birational
$G$-equivariant 
morphism $f: X_{\min}\to Y$ is an isomorphism.
In this situation $X_{\min}$ is called \emph{$G$-minimal model} of $X$.
If the surface $X$ is additionally rational, then 
 one of the following holds \cite{Iskovskikh-1979s-e}:
\begin{itemize}
 \item
$X_{\min}$ is a del Pezzo surface whose invariant Picard number $\Pic(X_{\min})^G$ is of rank $1$, or
\item
$X$ admits a structure of $G$-conic bundle, that is, there exists a
surjective $G$-equivariant morphism $f: X_{\min}\to \PP^1$ such that $f_*\OOO_{X_{\min}}=\OOO_{\PP^1}$,
$-K_{X_{\min}}$ is $f$-ample and $\rk \Pic(X_{\min})^G=2$.
\end{itemize}

\end{say}

\begin{say}{\bf Stable conjugacy.} 
We say that $G$-varieties
$(X,G)$ and $(Y,G)$ are \textit{stably birational}
if for some $n$ and $m$ there exists an equivariant 
birational map $X\times \PP^n \dashrightarrow Y\times \PP^m$, where 
actions on $\PP^n$ and $\PP^m$ are trivial.
This is equivalent to the conjugacy of subgroups
$G\subset \Bbbk(X)(t_1,\dots,t_n)$ and $G\subset \Bbbk(Y)(t_1,\dots,t_m)$.

By the \emph{no-name} lemma 
we have the following.

\begin{ppar}{\bf Remark.}\label{Remark-faithful}
Let $V$, $W$ be \emph{faithful 
linear} representations of $G$. 
Then the $G$-varieties $(V,G)$ and $(W,G)$
are stably conjugate. Indeed, let $n:=\dim V$, $m:=\dim W$. Consider trivial 
linear representations $V'$ and $W'$ with $\dim V'=n$ and $\dim W'=m$. 
According to the \textit{no-name lemma} (see e.g. \cite[Appendix 3]{Shafarevich1994a})
we can choose invariant coordinates for semi-linear action of $G$ on $V\otimes \Bbbk(W)$.
This means that two embeddings $G\subset \Cr_{n+m}(\Bbbk)$
induced by actions on $V\times W$ and $V'\times W$
are conjugate. Similarly,
the embeddings $G\subset \Cr_{n+m}(\Bbbk)$
induced by actions on $V\times W$ and $V\times W'$
are also conjugate. Hence, $(V,G)$ and $(W,G)$
are stably conjugate. 
\end{ppar} 

\begin{ppar}{\bf Definition.}\label{Definition-definitionslinearizable}
We say that a $G$-variety $(X,G)$
(or, by abuse of language, a group $G$) is \textit{stably linearizable}
if it is stably birational to $(V,G)$, where $V=\Bbbk^m$ is some faithful linear representation.
\end{ppar}

\begin{ppar}{\bf Remark.}\label{Remark-definitions}
One can define stable linearizability is several other ways:
\begin{enumerate}
\item\label{Remark-definitions-i}
if $(X,G)$ is stably birational to $(\PP^N,G)$ for some $N$;
\item\label{Remark-definitions-ii}
if $(X,G)$ is stably birational to $(\PP^N,G)$ for  $N=\dim X$;
\item\label{Remark-definitions-iii}
if there exists a $G$-birational map $X\times \PP^n \dashrightarrow  \PP^N$ for some $N$,
where the action on $\PP^n$ is trivial.
\end{enumerate}
In the view of Remark \ref{Remark-faithful} our definition \ref{Definition-definitionslinearizable}
seems to be most natural one. Clearly, we have the following implications:
\[
\text{\ref{Definition-definitionslinearizable}}
\Longrightarrow \text{\ref{Remark-definitions-iii}}\Longrightarrow
\text{\ref{Remark-definitions-i}},\quad \text{\ref{Remark-definitions-ii}}\Longrightarrow
\text{\ref{Remark-definitions-i}}
\]
The example below shows that, in general, the implications 
\ref{Remark-definitions-i}, \ref{Remark-definitions-ii}, \ref{Remark-definitions-iii}
$\Longrightarrow$
\ref{Definition-definitionslinearizable} do not hold.
\end{ppar}


\begin{ppar}{\bf Example.}\label{Example-one}
Let $\mathbf Q_8$ be the quaternion group of order $8$
and let $V$ be its faithful two-dimensional irreducible representation.
Then, for any $r$, the $(2r-1)$-dimensional projective space 
$\PP(V^{\oplus r})$ is a $G$-variety, where
$G=\mathbf Q_8/[\mathbf Q_8,\mathbf Q_8]\simeq (\ZZ/2\ZZ)^2$.
It is easy to see that there is no fixed points on this $\PP(V^{\oplus r})$.
Applying Lemma \ref{Lemma-Kollar-Szabo} (below) one can see that 
the $G$-variety $(\PP^{2r-1}, G)$ is not stably 
linearizable.
Similar examples can be constructed for the group $G=(\ZZ/n\ZZ)^2$ 
(e.g. instead of $\mathbf Q_8$ one can start with
the Heisenberg group of order $p^3$).
\end{ppar}

\begin{stheorem}{\bf Lemma (see \cite{Kollar-Szabo-2000}).}
\label{Lemma-Kollar-Szabo}
For any finite abelian group $G$ and any
$G$-birational map $X\dashrightarrow Y$ of complete 
$G$-varieties, the set $X^G$ is non-empty if and only if 
so $Y^G$ is. 
\end{stheorem}
\end{say}

\begin{say}{\bf Stable conjugacy and $H^1(G,\Pic(X))$.}

\begin{ppar} {\bf Definition.}
We say that a non-singular $G$-variety $(X,G)$ is \textit{$H^1$-trivial} if
$H^1(G_1,\Pic(X))=0$ for any subgroup $G_1\subset G$.
\end{ppar}

\begin{stheorem}{\bf Theorem (\cite{Bogomolov-Prokhorov}).}
Let $(X,G)$ be a smooth projective $G$-variety.
If $(X,G)$ is stably linearizable, then
$(X,G)$ is $H^1$-trivial.
\end{stheorem}
Note that the inverse implication does not true in general 
(see Remark \ref{Remark2-Iskovskikh-surface}).
Note also that the assertion of the theorem holds for 
any other definition of stable linearizability 
\ref{Remark-definitions}\ref{Remark-definitions-i}-\ref{Remark-definitions-iii}.

Our basic tool is the following theorem proved in \cite{Bogomolov-Prokhorov}.
\begin{stheorem}{\bf Theorem {\cite{Bogomolov-Prokhorov}}.}
\label{theorem-main-p}
Let $(X,G)$ be a non-singular projective rational $G$-surface, where 
$G$ is a cyclic group $G$ of prime order $p$.
Assume that $G$ fixes \textup(point-wise\textup) a curve of genus $g>0$.
Then
\begin{equation*}
H^1(G,\Pic(X))\simeq (\ZZ/p\ZZ)^{2g}.
\end{equation*}
If $H^1(G,\Pic(X))=0$, then $(X,G)$ is linearizable. 
\end{stheorem}

\begin{stheorem}{\bf Corollary.}\label{Corollary-H1}
Let $(X,G)$ be a non-singular projective rational $G$-surface, where 
$G$ is an arbitrary finite group.
If $(X,G)$ is $H^1$-trivial, then any non-trivial element of 
$G$ does not fix a curve of positive genus.
\end{stheorem}

\end{say}

\section{Group actions on del Pezzo surfaces}
\begin{say}
Let $X$ be a del Pezzo surface of degree $d\le 6$, i.e. $K_X^2=d$.
It is well-known that $X$ can be realized as the blowup $X\to \PP^2$
of $r:=9-d$ points in general position.
The group $\Pic(X)\simeq \ZZ^{r+1}$ has a basis $\mathbf h,\mathbf e_1,\dots, \mathbf e_r\in \Pic(X)$,
where $\mathbf h$ is the pull-back of the class of a line on $\PP^2$ and the $\mathbf e_i$'s
are the classes of exceptional curves.
\end{say}

\begin{say}
\label{delta}
Put
\[
\Delta_r:=\{\mathbf{x} \in \Pic(X) \mid \mathbf{x}^2=-2,\quad \mathbf{x}\cdot K_X=0\}.
\]
Then $\Delta_r$ is a root system in the orthogonal complement to $K_X$ in $\Pic(X)\otimes \RR$. 
Depending on $d$, the type of
$\Delta_r$ is the following (\cite{Manin-Cubic-forms-e-I}):
\bigskip
\begin{center}
\begin{tabular}[]{c|cccccc}
\hline
\\[-8pt]
$d$ &1&2&3&4&5&6
\\
\hline
\\[-5pt]
$\Delta_r$ & $\mathrm E_8$& $\mathrm E_7$& $\mathrm E_6$& $\mathrm D_5$& $\mathrm A_4$& $\mathrm A_1\times \mathrm A_2$
\\
\hline
\end{tabular}
\end{center}
\end{say}

\begin{ppar} \label{ppar-equation-injection}
There is a natural homomorphism 
\begin{equation}\label{equation-injection}
\varrho: \Aut(X) \longrightarrow \mW(\Delta_r),
\end{equation}
where $\mW(\Delta_r)$ is the Weyl group of $\Delta_r$. This homomorphism 
is injective if $d\le 5$
(see e.g. \cite[Corollary 8.2.32]{Dolgachev-topics}).
\end{ppar}

Denote by
$Q=\mQ(\Delta_r)$ the sublattice of $\Pic(X)$ generated by the roots.
Clearly, $\mQ(\Delta_r)$ coincides with the lattice of integral points
in $K_X^\perp\subset \Pic(X)\otimes \RR$.

\begin{say}\label{subsection-Lefschetz}
For an element $\delta\in \mW(\Delta_r)$ or $\Aut(X)$,
denote by $\tr(\delta)$ its trace on $Q$.
Let $G\subset \Aut(X)$ be a (finite) subgroup and let $n$ be the order of $G$. 
Computing the character of the trivial subrepresentation we get
\begin{equation}
\label{equation-Lefschetz-1}
\rk \Pic(X)^G= 1+\frac 1{n} \sum_{\delta\in G} \tr(\delta).
\end{equation}
On the other hand, since $\operatorname{Tr}_{H^2(X,\RR)}(\delta)=1+\tr(\delta)$,
by the Lefschetz fixed point formula we have
\begin{equation}
\label{equation-Lefschetz-2}
\Eu(X^\delta)=\tr(\delta)+3.
\end{equation}
\end{say}

\begin{say}
Now we prove the implication \ref{Theorem-main-del-Pezzo-3}$\Rightarrow$\ref{Theorem-main-del-Pezzo-1}
of Theorem \xref{Theorem-main-del-Pezzo}.
By \cite[Proposition 31.3]{Manin-Cubic-forms-e-I} we have the following.
\begin{stheorem}{\bf Corollary.}\label{Corollary-Manin-computation}
Let $(X,G)$ be a projective $G$-surface.
Let $\{ C_i\}$ be a finite $G$-invariant set of irreducible
curves whose classes generate $\Pic(X)$.
If $G$ acts on $\{ C_i\}$ transitively,
then $H^1(G,\Pic(X))=0$.
\end{stheorem}

\begin{stheorem}{\bf Proposition.}\label{lemma-H1-trivial-d>5}
Let $(X,G)$ be a projective non-singular rational surface with $K_X^2\ge 5$.
Then $H^1(G,\Pic(X))=0$.
\end{stheorem}

\begin{proof}
To show that $H^1(G,\Pic(X))=0$ we may assume that $(X,G)$ is $G$-minimal
(otherwise we replace $X$ with its minimal model).
If $K_X^2\ge 8$, then $X$ is either $\PP^2$ or a Hirzebruch surface $\FF_e$
and $G$ acts on $\Pic(X)$
by (possibly trivial) permutation of the extremal rays.
Hence, $\Pic(X)$ is a permutation $G$-module and 
$H^1(G,\Pic(X))=0$. Thus $K_X^2=6$ or $5$ and $X$ is a del Pezzo surface with
$\rk \Pic(X)^G=1$ (see \cite{Iskovskikh-1979s-e}).

If $K_X^2=6$, then $X$ contains exactly $6$ lines $C_1,\dots,C_6\subset X$.
Since $\Pic(X)^G=\ZZ\cdot K_X$, these lines form one $G$-orbit.
By Corollary \ref{Corollary-Manin-computation} we conclude that $H^1(G,\Pic(X))=0$.

Finally, consider the case $K_X^2=5$.
Then $\Aut(X)\simeq \mW(\mathrm A_4)\simeq \mathfrak S_5$
(see e.g. \cite[Theorem 8.5.8]{Dolgachev-topics}).
Let $\LLL:=\{L_1,\dots,L_{10}\}$ be the set of lines on $X$.
The action of $G$ on $\LLL$ is faithful (see \ref{ppar-equation-injection}).
Let $\LLL=O_1\cup\cdots\cup O_l$ be the decomposition in
$G$-orbits and let $r_i$ be the cardinality of $O_i$.
Then $\sum r_i=10$.
Since $\Pic(X)^G=\ZZ\cdot K_X$, each number $r_i$ is divisible by $5$.
By Corollary \ref{Corollary-Manin-computation} we have only one possibility:
$r_1=r_2=5$. In particular, the order of $G$ is divisible by $5$.
Then both $O_1$ and $O_2$ form anticanonical divisors and 
the corresponding dual graphs
are combinatorial cycles.
In this case $G$ contains no elements of order $3$.
Hence the order of $G$ divides $20$ and $G$ has a normal subgroup $\langle\delta\rangle$ of 
order $5$. Since $\tr(\delta)=-1$, by the Lefschetz fixed point formula 
$\Eu(X^\delta)=2$. Write $X^\delta=V_1\cup V_0$, where $V_0\cap V_1=\emptyset$,  $\dim V_0=0$,
and $V_1$ is of pure dimension one. 
The action of $G$ preserves this decomposition.
If $V_1\neq\emptyset$, then $V_1$ meets the cycle of lines corresponding to 
$O_1$. But then $\delta$ acts on $O_1$ trivially, a contradiction.
Hence, $V_1\neq\emptyset$ and so 
$\delta$ has exactly two isolated fixed points $P_1,\, P_2\in X$.
By blowing $\{P_1,P_2\}$ up we get a cubic surface $\tilde X$ containing $G$-invariant pair of skew lines.
Then 
a well-known classical construction gives us a birational equivariant transformation 
$\tilde X\dashrightarrow \PP^1\times \PP^1$
(cf. \cite[\S 8]{Dolgachev-Iskovskikh}).
Then by the above considered case $K_X^2=8$ we have 
$H^1(G,\Pic(X))=0$.
\end{proof}
\begin{stheorem}{\bf Corollary.}\label{Corollary-lemma-H1-trivial-d=4}
 Let $(X,G)$ be a $G$-del Pezzo surface described in \eqref{equation-d=4}-\eqref{equation-d=4-action}.
Then 
$(X,G)$ is $H^1$-trivial.
\end{stheorem}
\begin{proof}
If $G'\subset G$ is a proper subgroup, then 
$(X,G')$ is not minimal and
$H^1(G',\Pic(X))=0$ by Proposition \ref{lemma-H1-trivial-d>5}.
It is easy to see that the set of lines on $X$ has exactly two $G$-orbits
consisting of $4$ and $12$ elements.
Then $H^1(G,\Pic(X))=0$ by 
\cite[Ch. 4, \S 31, Table 2]{Manin-Cubic-forms-e-I}.
\end{proof}
\end{say}
\begin{say}\label{assumption-main-del-Pezzo}
The implication \ref{Theorem-main-del-Pezzo-2}$\Rightarrow$\ref{Theorem-main-del-Pezzo-3}
of Theorem \xref{Theorem-main-del-Pezzo}
is an immediate consequence of the following proposition which will be proved below in \S\ref{section-d=4}-\S\ref{section-d=1}.
\begin{stheorem}{\bf Proposition.}\label{Proposition-main-del-Pezzo}
Let $(X,G)$ is a minimal $G$-del Pezzo surface of degree $\le 4$
such that any non-identity element of $G$ does not fix a curve of positive genus.
Then $(X,G)$ is isomorphic to a $G$-surface described in \eqref{equation-d=4}-\eqref{equation-d=4-action}.
\end{stheorem}
\end{say}

\section{Quartic del Pezzo surfaces}\label{section-d=4}
\begin{say}
Throughout this section, let $X$ be a del Pezzo surface of degree $4$.
It is well-known that the anti-canonical linear system
embeds $X$ to $\PP^4$ so that the image is a complete
intersection of two quadrics.
In a suitable coordinate system in $\PP^4$ the equations of $X$ can be written
in the form
\begin{equation}
\label{equation-2-quadrics}
\sum_{i=0}^4 x_i^2= \sum_{i=0}^4 \theta_i x_i^2 =0,
\end{equation}
where the $\theta_i$'s are distinct constants
(see e.g. \cite[Lemma 8.6.1]{Dolgachev-topics}).
We regard these constants $\theta_i\in \Bbbk$ as points of a projective line.
In other words, quadrics passing through $X$ form a pencil 
$\mathscr Q$ and the points $\theta_i$ correspond 
to degenerate members of $\mathscr Q$.
Five commuting involutions
$\tau_i: x_i \mapsto -x_i$ generate a normal abelian subgroup
$A\subset \Aut(X)$ with a unique relation
$\tau_1\cdots \tau_5=\id$. 
Thus
\[
A=\{1,\, \tau_k,\, \tau_i\tau_j \mid 1\le k\le 5,\ 1\le i < j\le 5\},
\qquad A\simeq (\ZZ/2\ZZ)^4.
\]
\end{say}

\begin{say}\label{Weyl-group-D5}
It is well-known (see e.g. \cite{Bourbaki2002}) that the root system of type $\mathrm D_5$
can be realized as the set $\pm \mathbf r_i \pm \mathbf r_j$,
where $\mathbf r_1,\dots,\mathbf r_5$ is the standard basis of
$\mathbb R^5$. The Weyl group
$\mW(\mathrm D_5)$ is the semi-direct product
$(\ZZ/2\ZZ)^4 \rtimes \Sym_5$, where $(\ZZ/2\ZZ)^4$ acts on $\mathbb R^5$
by $\mathbf r_i \mapsto (\pm 1)_i \mathbf r_i$ so that $\prod_i (\pm 1)_i=1$
and $\Sym_5$ acts on $\mathbb R^5$
by permutations of the $\mathbf r_i$'s.

The image $\varrho(A) \subset \mW(\mathrm D_5)$
under the injection
\eqref{equation-injection}
coincides with
$(\ZZ/2\ZZ)^4 \subset (\ZZ/2\ZZ)^4\rtimes \Sym_5$.
Thus we identify
$\varrho(A)$ with $(\ZZ/2\ZZ)^4$
and $\varrho(\tau_i)$ with $\tau_i$.
Note the fixed point locus of each $\tau_i$ is an elliptic curve
that cut out on $X$ by the hyperplane $\{x_i=0\}$
(and so the $\tau_i$'s are de Jonqui\`eres involutions of genus $1$).
The fixed point loci of other involutions in $A$
consist of exactly four points. Therefore,
\begin{equation}
\label{equation-deg=4-tr} 
\tr(\tau_{i})= -3\ \forall i,\qquad \tr(\tau_{i}\tau_{j})= 1\ \forall i\neq j.
\end{equation}
\end{say}

\begin{say}\label{subsection-d=4-pencils}
Another, intrinsic description of the $\tau_i$'s
is as follows. On $X$ there are 10 pencils of conics
$\mathscr C_1,\dots,\mathscr C_5,\mathscr C_1',\dots,\mathscr C_5'$
so that these pencils satisfy the conditions $\mathscr C_i\cdot \mathscr C_i'=2$,\
$\mathscr C_i\cdot\mathscr C_j=\mathscr C_i\cdot\mathscr C_j'=1$ for $i\neq j$,
and $\mathscr C_i+ \mathscr C_i'\sim -K_X$.
Two ``conjugate'' pencils $\mathscr C_i$ and $\mathscr C_i'$ define
a double cover $\psi_i: X\to \PP^1\times \PP^1$. Then $\tau_i$ is
the Galois involution of $\psi_i$. Note that $\psi_i$ coincides with
the projection of $X$ from the vertex of
a singular quadric of the pencil generated by \eqref{equation-2-quadrics}.
Thus there are the following canonical bijections:
\begin{equation}
\label{equation-d=4-tau}
\{\tau_i\} \longleftrightarrow \{\psi_i\} \longleftrightarrow
\{ (\mathscr C_i,\mathscr C_i')\}\longleftrightarrow \{\theta_i\},\quad i=1,\dots,5.
\end{equation}
The group $\Aut(X)$ acts on the pencil
of quadrics $\mathscr Q_\lambda$ in $\PP^4$ generated by \eqref{equation-2-quadrics}
so that the set of degenerate quadrics corresponding to the values
$\lambda=\theta_i$, $i=1,\dots,5$ is preserved.
Hence there exist homomorphisms
\[
\varrho_1:\Aut(X)\to PGL_2(\Bbbk),\qquad \varrho_2:\Aut(X)\to \Sym_5
\]
with $\ker(\varrho_1)=\ker(\varrho_2)=A$. 
This immediately gives us
the following possibilities for the group $\Aut(X)/A$ (see \cite[\S 6]{Dolgachev-Iskovskikh}):
\begin{equation}\label{equation-d=4-possibilities}
\{1\},\ \ZZ/2\ZZ,\ \ZZ/3\ZZ,\ \ZZ/4\ZZ,\ \ZZ/5\ZZ,\ \Sym_3,\ \DD_{5}.
\end{equation}
\end{say}

\begin{say}\label{notation-d=4-r_i}
Now let a finite group $G$ faithfully act on $X$ so that $(X,G)$ is minimal (i.e. $\Pic(X)^G\simeq \ZZ$)
and any non-identity element of $G$ does not fix a curve of positive genus.
Denote $A_G:=G\cap A$.
For short, we identify $\varrho(G)$ with $G$.

Recall that $K_X^2=4$.
Let $\LLL:=\{L_1,\dots,L_{16}\}$ be the set of lines on $X$.
Let $\LLL=O_1\cup\cdots\cup O_l$ be the decomposition in
$G$-orbits and let $r_i$ be the cardinality of $O_i$.
Then $\sum r_i=16$.
Since $\Pic(X)^G=\ZZ\cdot K_X$, each number $r_i$ is divisible by $4$.
\end{say}
By our assumption in \ref{notation-d=4-r_i} we have the following.
\begin{stheorem}{\bf Corollary.}
$G\not\ni \tau_i$ for $i=1,\dots,5$.
\end{stheorem}

The following lemma is an immediate consequence of the description of $A$.
\begin{stheorem}{\bf Lemma.}\label{Corollary-d=4-A}
There are two kinds of non-trivial subgroups $A'\subset A$
satisfying the property $A'\not\ni \tau_i$ for $i=1,\dots,5$:
\begin{itemize}
 \item
 $A_{i,j}=\{1,\, \tau_i\tau_j \mid i\neq j\}$, and
 \item
$A_{k,l,m}=\{1,\, \tau_k\tau_l,\, \tau_l\tau_m,\, \tau_k\tau_m \mid k\neq l\neq m\neq k\}$.
\end{itemize}
\end{stheorem}
\begin{ppar}{\bf Remark.}
Note that if $A_G=A_{i,j}$, then $A_G$ is contained in the center of $G$.
Using \eqref {equation-deg=4-tr} we immediately conclude that
\begin{equation}
\label{equation-d=4-tr-A}
\sum_{\upsilon\in A_G} \tr(\upsilon)=
\begin{cases}
6&\text{if $A_G=A_{i,j}$},
\\
8&\text{if $A_G=A_{k,l,m}$}.
\end{cases}
\end{equation}
\end{ppar}

For $G/A_G$ we have the same possibilities \eqref{equation-d=4-possibilities}
as for $\Aut(X)/A$. 
Consider these possibilities case by case. By \eqref{equation-d=4-tr-A}
and \eqref{equation-Lefschetz-1} \ $G\neq A_G$.
\begin{say}
{\bf Cases $G/A_G \simeq \ZZ/5\ZZ$ and $\DD_{5}$.} 
The order of $G$ divides $40$. By Sylow's theorem
the Sylow $5$-subgroup $G_5\subset G$ is normal.
By \ref{notation-d=4-r_i} we see that $r_i\not\equiv 0\mod 5$ for all $i$.
Hence 
$G_5$ is contained in the stabilizer of any line $L\in \LLL$.
But then the action of $G$ on $\LLL$ and on $\Pic(X)$ is not faithful, a
contradiction.
\end{say}

\begin{say} {\bf Case $G/A_G \simeq \ZZ/3\ZZ$.}
For convenience of the reader we reproduce here the following fact from \cite[\S 6]{Dolgachev-Iskovskikh}:
\begin{stheorem}{\bf Lemma \cite[\S 6]{Dolgachev-Iskovskikh}.}\label{lemma-description-d=4-S3}
Let $X$ be a quartic del Pezzo surface and let $\gamma \in \Aut(X)$ be an element of order $3$.
Then $X$ is isomorphic to the surface given by \eqref{equation-d=4}.
Moreover, $\Aut(X)\simeq A\rtimes \Sym_3$.
The center of $\Aut(X)$ is of order $2$ and generated by an element of the form $\tau_i\tau_j$, $i\neq j$.
\end{stheorem}

\begin{proof}
Since $X$ contains exactly $16$ lines, there exists 
at least one $\gamma$-invariant line $L\subset X$.
Let $L_1,\dots,L_5\subset X$ be (skew) lines meeting $L$ and
let $f:X\to \PP^2$ be the contraction of $L_1,\dots,L_5$.
Let $C:=f(L)$ and $P_i=f(L_i)$. Then the action of $\gamma$ on $X$ 
is induced by one on $C\subset \PP^2$. Up to permutation 
of $L_1,\dots,L_5$ we may assume that $\gamma$ fixes $P_1$ and $P_2$
and permutes $P_3,\, P_4,\, P_5$. Then the set $\{P_1,\dots,P_5\}$ 
is unique up to projective equivalence. Hence $X$ is unique up to
isomorphism. On the other hand, it is easy to see that 
the surface \eqref{equation-d=4} admits an isomorphism $\gamma$ of order $3$
given by \eqref{equation-d=4-action}.
Moreover, $\Aut(X)$ contains the group $A\rtimes \Sym_3$
generated by $A$, $\gamma$ and 
\[
\beta: (x_1,x_2,x_3,x_4, x_5) \longmapsto (x_1, x_3,x_2,x_5,x_4).
\]
By \eqref{equation-d=4-possibilities} we see that $\Aut(X)=A\rtimes \Sym_3$.
\end{proof}

\begin{stheorem}{\bf Corollary.}\label{Claim-deg=4-order=3}
Let $\gamma\in \Aut(X)$ be an element of order $3$.
Then $X^\gamma$ consists of exactly $5$ points. 
\end{stheorem}

By Corollary \ref{Claim-deg=4-order=3} the exists a 
$G$-fixed point $P\in X$. Since in a neighborhood of $P$ 
the action of $(\ZZ/2\ZZ)^2$ cannot be free in codimension one,
we have $A_G=A_{i,j}$ for some $i\neq j$. 
Hence $G$ is cyclic of order $6$. Since the cardinality of any 
orbit $O_i\subset \LLL$ must be divisible by $4$, we get a contradiction.
\end{say}

\begin{say}\label{d=4-case-realized}
{\bf Case $G/A_G \simeq \Sym_3$. }
We show that only the possibility \ref{Theorem-main-del-Pezzo-3}(b) of Theorem \ref{Theorem-main-del-Pezzo} occurs here.
Let $G_3$ (resp. $G_2$) be a Sylow $3$ (resp. $2$)-subgroup
of $G$. Clearly, $G_2\supset A_G$ and $G_2/A_G\simeq \ZZ/2\ZZ$.
By Lemma
\ref{lemma-description-d=4-S3} \ $X$ is isomorphic to the surface given by \eqref{equation-d=4}, 
$\Aut(X)\simeq A\rtimes \Sym_3$, and the center of $\Aut(X)$ is 
generated by an element $\tau_i\tau_j$, $i\neq j$.

\begin{stheorem}{\bf Lemma.}\label{lemma-description-d=4-S3-aut}
In the above settings the image of the natural representation 
$\varrho: \Aut(X)\hookrightarrow \mW(\mathrm D_5)\subset GL(Q)$ is 
contained in $SL(Q)$.
\end{stheorem}
\begin{proof}
 By \ref{Weyl-group-D5} we can write the elements of $A$ in a diagonal form 
so that $A\subset SL(Q)$ and the determinant of any element of 
$\mW(\mathrm D_5)$ equals to $\pm 1$.
The fixed point locus of $\beta$ consists of a smooth rational curve 
and a pair of isolated points. Hence,
$\tr(\beta)=1$ and so $\det(\beta)=1$.
This implies that the image of the whole group $\Aut(X)$ 
is contained in $SL(Q)$. 
\end{proof}

\begin{ppar}
Assume that $A_G=A_{i,j,k}$. Since elements of $A_G$ and $G_3$ do not commute,
$G_3$ is not normal in $G$.
By Sylow's theorem the number
of Sylow $3$-subgroups equals to $4$.
The action on the set of these subgroups induces an isomorphism
$G\simeq \Sym_4$. 
By Corollary \ref{Claim-deg=4-order=3} for the elements $\gamma\in G$ of order $3$ we have
$\tr(\gamma)=2$. Hence, by \eqref{equation-d=4-tr-A} and \eqref{equation-Lefschetz-1}
\[
 \sum_{\upsilon\in \Alt_4} \tr(\upsilon)= 24,\qquad \sum_{\upsilon\in \Sym_4\setminus \Alt_4} \tr(\upsilon)= -24.
\]
Since $\Eu(X^\upsilon)>0$ for all $\upsilon\in G$, we have $\tr(\upsilon)=-2$ 
for all $\upsilon\in \Sym_4\setminus\Alt_4$.
In our case, $\dim Q=5$. 
Hence $\tr(\upsilon)$ must be odd for an element of order $2$, a contradiction. 
\end{ppar}

\begin{ppar}\label{d=4-S3-Aut}
Thus $A_G=A_{i,j}$. Then 
$G_3$ is normal in $G$ and so 
$G$ is a semi-direct product $G=G_3\rtimes G_2$ which is not a direct product because 
$G$ is not abelian. For short, we identify $G$ with its image in $\mW(\mathrm D_5)\subset GL(Q)$.
We claim that $G_2$ is cyclic.
Indeed, otherwise $G\simeq \Sym_3\times (\ZZ/2\ZZ)$.
It is easy to check that in this case $Q$ must contain a trivial $G$-representation 
(because $G\subset SL(Q)$ by Lemma \ref{lemma-description-d=4-S3-aut}). Since $\Pic(X)^G\simeq \ZZ$,
this is impossible. Therefore,
$G_2\simeq \ZZ/4\ZZ$ and $G\simeq (\ZZ/3\ZZ) \rtimes(\ZZ/4\ZZ)$.
Up to permutations of coordinates we may assume that the center of $\Aut(X)$
is generated by 
\[
\delta=\tau_4\tau_5: (x_1,x_2,x_3,x_4, x_5) \longmapsto(x_1,x_2,x_3,-x_4, -x_5).
\]
Clearly, the center of $G$ commutes with all elements of $\Aut(X)$. 
Thus $\delta\in G$. 
Now let $\beta^\bullet$ (resp. $\gamma^\bullet$) be an element of 
$G$ of order $4$ (resp. $3$) whose image in $\Sym_3$ coincides with 
$\beta$ (resp. $\gamma$). Thus $\beta^\bullet(x_i)=\pm \beta(x_i)$ and 
$\gamma^\bullet(x_i)=\pm \gamma(x_i)$ for all $i$. 
Since ${\gamma^\bullet}^3=\id$, replacing $x_i$ with $\pm x_i$ 
we may assume that 
$\gamma^\bullet =\gamma$.
Since $(\beta^\bullet)^2=\delta$ 
and $\beta^\bullet\gamma{\beta^\bullet}^{-1}=\gamma^{-1}$,
as above, 
we get $\beta^\bullet=\beta'$.
Thus our group $G$ coincides with that constructed 
in \eqref{equation-d=4}-\eqref{equation-d=4-action}.
It remains to show that this group is minimal.
Let $\nu\in G$ be an element of even order $2k$.
Then $\nu^k=\delta$ and so $X^\nu=(X^\delta)^\nu$.
Recall that $X^\delta$ is a set of four points.
Then one can easily see that $\Eu(X^\nu)=1$ (resp. $2$)
if $k=3$ (resp. $2$).
Thus we have 
\[
\sum_{\upsilon\in G} \tr(\upsilon)=
5+1+2\cdot 2- 2\cdot 2-6\cdot 1=0.
\]
By \eqref{equation-Lefschetz-1} we have
$\rk \Pic (X)^G =1$, i.e. $G$ is minimal.
\end{ppar}

\begin{ppar}{\bf Remark.}\label{Remark-d=4-case-realized}
Note that our group 
$G$ acts on $X^{G_3}$ and by Corollary \ref{Claim-deg=4-order=3}
there is a $G$-fixed point $P\in X^{G_3}$
so that $P$ does not lie 
on any line. Let $\tilde X\to X$ be the blowup of $P$.
Then $\tilde X$ is a cubic surface admitting an action of $G$ so that
$\rk \Pic(\tilde X)^G=2$. The exceptional divisor is an invariant line $L\subset\tilde X$
and the projection from $L$ gives a structure of $G$-equivariant conic bundle
$\tilde X\to \PP^1$. Thus we are in the situation described below in 
Theorem \ref{theorem-conic-bundle-1}
and Construction \ref{Construction-conic-bundle-non-exceptional} (with $n=3$).
\end{ppar}
\end{say}

\begin{say} {\bf Case $G/A_G \simeq \ZZ/2\ZZ$.}
Since $\Pic(X)^G\simeq \ZZ$, $A_G\neq \{1\}$.
Assume that $A_G=A_{i,j}$ for some $i$, $j$. Then by \eqref{equation-d=4-tr-A} we have 
$\sum_{\delta\in G\setminus A_G} \tr(\delta)=-6$.
Hence there exists $\delta\in G\setminus A_G$ such that 
$\Eu(X^{\delta})\le 0$. Since $X^{\delta}\neq\emptyset$, the element
$\delta$ fixes point-wise a curve of positive genus.
This contradicts our assumption \ref{notation-d=4-r_i}.
Therefore, $A_G=A_{i,j,k}$ for some $i$, $j$, $k$. 
In particular, $G$ is a (non-cyclic) group of order $8$.
Again by \eqref{equation-d=4-tr-A} we have 
$\sum_{\delta\in G\setminus A_G} \tr(\delta)=-8$
and $\Eu(X^{\delta})> 0$ for all $\delta\in G\setminus A_G$.
Hence, $\Eu(X^{\delta})=1$ for all $\delta\in G\setminus A_G$.
This means that any element $\delta\in G\setminus A_G$ has a
unique fixed point and the action of $G$ on $X$ is free in codimension $1$.
Applying the holomorphic Lefschetz fixed point formula,
we obtain that any $\delta\in G\setminus A_G$ has at least two
fixed points, a contradiction. 
\end{say}

\begin{say} {\bf Case $G/A_G \simeq \ZZ/4\ZZ$.}
Note that the stabilizer of $A_{i,j}$ (and $A_{k,l,m}$)
in $\Sym_5=\mW(\mathrm D_5)/A$ is the group $\Sym_2\times \Sym_3$.
Hence neither $A_{i,j}$ nor $A_{k,l,m}$ can be a normal subgroup of $G$.
Thus, $A_G=\{1\}$. Again we have
$0=5+\tr(\delta^2)+2\tr(\delta)$, where $\tr(\delta)$, $\tr(\delta^2)\ge -2$ by 
\eqref{equation-Lefschetz-2}
because $G$ does not fix a curve of positive genus.
We get only one possibility: $\tr(\delta^2)=-1$, $\tr(\delta)=-2$.
Hence, $X^G$ is a point, say $P$, and $X^{\delta^2}$ is either a smooth rational curve
or a pair of points. On the other hand, $X^{\delta^2}\ni P$
and $G$ acts on $X^{\delta^2}$ fixing $P$, a contradiction.
\end{say}
Thus Proposition \ref{Proposition-main-del-Pezzo}
is proved in the case $K_X^2=4$.

\section{Cubic surfaces}\label{section-cubics}
\begin{say}
Throughout this section $X$ denotes a cubic surface $X\subset \PP^3$. 
Let $G\subset \Aut(X)$ be a subgroup such that $(X,G)$ is minimal and
any non-identity element of $G$ does not fix a curve of positive genus.
Since the embedding $X\subset \PP^3$ is anti-canonical, it is
$G$-equivariant.
By our assumption for any element $1\neq \delta\in G$
the set $(\PP^3)^{\delta}$ does not contain any hyperplane.
Let $\psi(x_1,x_2,x_3,x_4)=0$ be the equation of $X$.
We  choose homogeneous coordinates in $\PP^3$ so that $\delta$ has
a diagonal form. 
\end{say}

\begin{mtheorem}{\bf Claim.}\label{Claim-cubic-ord=2}
Let $\tau\in G$ be an element of order $2$.
Then in a suitable coordinates its action on $\PP^3$ has the form
$\tau=\diag(1,1,-1,-1)$ and 
\[
\psi=\psi_3(x_1,x_2)+ x_1\psi_2(x_3,x_4)+x_2\psi_2'(x_3,x_4),
\]
where $\deg \psi_3=3$, $\deg \psi_2=\deg \psi_2'=2$, and $\psi_3$ has no multiple factors.
Furthermore, $X^\tau=L(\tau) \cup \{P_1,\, P_2,\, P_3\}$,
where $L(\tau):=\{x_1=x_2=0\}$ and $\{P_1,\, P_2,\, P_3\}=X\cap \{x_3=x_4=0\}$. In particular,
$\Eu(X^\tau)=5$.
\end{mtheorem}
\begin{proof}
Since $(\PP^3)^{\tau}$ does not contain any hyperplane, we can write 
$\tau=\diag(1,1,-1,-1)$. 
Replacing $\tau$ with $-\tau$ we may assume that $\psi$ is 
invariant. The rest is obvious. 
\end{proof}
\begin{mtheorem}{\bf Claim.}\label{Claim-cubic-ord=3}
Let $\tau\in G$ be an element of order $3$.
Then the fixed point locus $X^\tau$ is zero-dimensional and $\Eu(X^\tau)\ge 3$.
\end{mtheorem}
\begin{proof}
Up to permutations of coordinates we may assume that $\delta$ has the form
$\diag(1,1,\zeta_3,\zeta_3)$ or $\diag(1,1,\zeta_3,\zeta_3^{-1})$.
Assume that $\dim X^\tau=1$. By the above there exists a line 
$L\subset X^\tau$.
It is well-known that a given line $L$ on a cubic surface meets 
exactly 10 other lines $L_1,\dots, L_{10}$ and up to reenumeration
one can assume that  the lines
 $\{L_1,\dots,L_5\}$ (resp. $\{L_6,\dots,L_{10}\}$)
are mutually disjoint. Then  each line $L_i$ must be 
$\delta$-invariant (because $L_i\cap L$ is a fixed point). 
In this case the classes of $L_1,\dots,L_5$ are contained in $\Pic(X)^\delta$ 
and linearly independent there. Since the canonical class $K_X$
is also $\delta$-invariant, we see that the action of $\delta$ on $\Pic(X)$
must be trivial. This contradicts the injectivity of  $\varrho: \Aut(X) \longrightarrow \mW(\operatorname{E}_6)$
(see \ref{ppar-equation-injection}).


Thus $\dim X^\tau=0$. On the other hand, $X^\tau\neq \emptyset$
and $\tr(\tau)=3$, $0$, or $-3$.
Hence, $\Eu(X^\delta)=6$ or $3$. 
\end{proof}

\begin{mtheorem}{\bf Lemma.}\label{Lemma-cubic-tr}
For any element $\delta\in G$ we have $\tr(\delta)\ge 0$
except for the following case:
\begin{enumerate}
\item[\rm (*)] 
$\ord(\delta)=6$, $\tr(\delta)=-1$, $X^\delta$ consists of two points: 
$X^\delta=L(\delta^3)^\delta=\{R_1,R_2\}$, where $L(\delta^3)$ is the line
introduced in Claim \xref{Claim-cubic-ord=2}.
Moreover, in the local coordinates near $R_i$
the action of $\delta^2$ is given by a scalar matrix.
\end{enumerate}
\end{mtheorem}

\begin{proof}
 By \cite{atlas} the orders of elements of $\mW(\mathrm E_6)$ 
are as follows:
$1$, $2$, $3$, $4$, $5$, $6$, $8$, $9$, $10$, $12$.
Consider the possibilities for $\delta\in G$.
Let $\chi(t)$ be the characteristic polynomial of $\delta$ on $Q$.
Clearly, $\deg \chi=6$ and $\chi$ is a product of cyclotomic polynomials $\Phi_d$, 
where $d$ divides $\ord(\delta)$.

If $\ord(\delta)\le 3$, then $\tr(\delta)\ge 0$ by Claims \ref{Claim-cubic-ord=2}
and \ref{Claim-cubic-ord=3}. Thus we may assume that $\ord(\delta)\ge 4$.
If $\ord(\delta)=5$, then 
the only possibility is $\chi=\Phi_{5}\Phi_1^2=t^6 - t^5 - t + 1$ and $\tr(\delta)=1$. 
If $\ord(\delta)=9$, then again we have 
$\chi=\Phi_9= t^6 + t^3 + 1$ and $\tr(\delta)=0$.

It remains to consider the case where the order of $\delta$ is even, so
$\ord(\delta)=2m$, $m=2$, $3$, $4$, $5$, or $6$.
Then $\delta^m$ is described in Claim \ref{Claim-cubic-ord=2} and so
\[
X^\delta=L^\delta \cup \{P_1,\, P_2,\, P_3\}^\delta.
\]
where $L:=L(\delta^m)$ and the points $P_1$, $P_2$, $P_3$ lie on one line in $\PP^3$.
Here $L^\delta$ either is a couple of points or coincides with $L$.
Hence, $\Eu(L^\delta)=2$ and $\{P_1,\, P_2,\, P_3\}^\delta=\emptyset$ if and only if 
$\delta$ permutes the $P_i$'s. Thus $\Eu(X^\delta)\le 2$ only if 
$m=3$, $\tr(\delta)=-1$, and $X^\delta=L^\delta$.
Consider the blow-down $X\to X'$ of $L$ to a point, say $R$. 
Since $\delta^2$ acts on $X$ freely in codimension one (see Claim \ref{Claim-cubic-ord=3}),
in the local coordinates near $R$
the action of $\delta^2$ can be written as $\diag(\zeta_3, \zeta_3^{-1})$.
Then it is easy to see that in the local coordinates near $R_i$
the action can be written as $\diag(\zeta_3^{k}, \zeta_3^{k})$, $k=1$ or $2$.
\end{proof}
\begin{say}\textit {Proof of Proposition \xref{Proposition-main-del-Pezzo} 
\label{d=3-finish}
in the case $K_X^2=3$.}
Since $(X,G)$ is minimal, we have $\sum_{\delta\in G} \tr(\delta)=0$ by
\eqref{equation-Lefschetz-1}.
Hence, $\tr(\delta)< 0$ for some $\delta\in G$.
By Lemma \ref{Lemma-cubic-tr} we have $\ord(\delta)=6$ and $\tr(\delta)=-1$.
Let $G_1,\dots, G_r\subset G$ be all cyclic subgroups generated by
such elements $\delta_i$ of order $6$.
We claim that $\delta_i^3\neq \delta_j^3$ for $i\neq j$.
Assume the converse: $\delta_i^3= \delta_j^3:=\tau$.
The element $\tau$ is described in Claim \ref{Claim-cubic-ord=2}.
Put $L:=L(\tau)$.
The projection from $L$ defines a $\langle \delta_i,\delta_j\rangle$-equivariant 
conic bundle structure $f: X\to\PP^1$
so that the restriction $f|_{L}: L\to \PP^1$ is a double cover.
It has two ramification points $R_1$, $R_2\in L$. Since each $\delta_i$ has 
exactly two fixed points, we have $X^{G_i}=X^{G_j}=\{R_1,\, R_2\}$.

Replacing $\delta_j$ with $\delta_j^{\pm 1}$ 
we may assume that the action of 
$\delta_i^2$ and $\delta_j^2$ on $T_{R_1,X}$ has the form 
$\diag(\zeta_3, \zeta_3)$. Hence, $\delta_i^2=\delta_j^2$
and so $\delta_i=\delta_j$ which proves our claim.
In particular, we see that for $i\neq j$ the intersection $G_i\cap G_j$ does not contain
any elements of order $2$.
Then by \eqref{equation-Lefschetz-1}
\[
0=\sum_{\delta\in G}\tr(\delta)>
\sum_{i=1}^r \left( \tr(\delta_i)+\tr(\delta_i^{-1})+\tr(\delta^3)\right)= 0.
\]
The contradiction proves 
Proposition \ref{Proposition-main-del-Pezzo} in the case $K_X^2=3$.
\hfill
$\square$
\end{say}

\section{Del Pezzo surfaces of degree $2$}\label{section-d=2}
\begin{say}
Throughout this section $X$ denotes a del Pezzo surface of degree $2$.
Recall that the anti-canonical map
is a double cover $X\to \PP^2$ branched over a smooth quartic $R\subset \PP^2$.
Let $\psi(x_0,x_1,x_2)=0$ be the equation of $R$.
Then $X$ can be given by the equation $y^2=\psi(x_0,x_1,x_2)$
in the weighted projective space $\PP(1,1,1,2)$.
The Galois involution $\gamma: X\to X$ of the double cover $X\to \PP^2$ is called
the \textit{Geiser
involution}. It is contained in the center of $\Aut(X)$
and $X^\gamma$ is a curve of genus $3$.
For any $\mathbf x\in \Pic(X)$ the
element $\mathbf x +\gamma^*\mathbf x$ is the pull-back of some
element of $\Pic(\PP^2)$. 
\end{say}

By \eqref {equation-Lefschetz-1} (cf. \ref{d=3-finish}) to establish 
Proposition \ref{Proposition-main-del-Pezzo} in the case $K_X^2=2$
it is sufficient to prove 
the following.
\begin{mtheorem}{\bf Lemma.}\label{lemma-d=2}
Let $G\subset \Aut(X)$ be a finite subgroup
such that any non-identity element of $G$ does not fix a curve of positive genus.
Then $\tr(\delta)\ge 0$ for any $\delta\in G$.
\end{mtheorem}

\begin{proof}
It is known that the center of $\mW(\mathrm E_7)$ is a cyclic group of order $2$
generated by the element $\gamma$ which 
is induced by the Geiser involution of $X$ and that acts as minus identity on $\mQ(\mathrm E_7)$. 
The quotient
$\mW(\mathrm E_7)/\langle\gamma \rangle$ is the (unique) simple group of order $1451520$
isomorphic to $PSp_6(\mathbb F_2)$. Let 
$\bar G$ be the image of $G$ in $\mW(\mathrm E_7)/\langle\gamma \rangle$.
By our assumption the group
$G$ does not contain $\gamma$. Hence, $G\simeq \bar G$.
Using the description of conjugacy classes in $PSp_6(\mathbb F_2)$
(see \cite{atlas}) we obtain that order of any element of $G$ is one of the following numbers:
1, 2, 3, 4, 5, 6, 7, 8, 9, 10, 12, 15.
Consider these possibilities case by case.
Let $\chi_\delta(t)$ denote the characteristic polynomial of the action of $\delta\in G$ on
$Q\otimes \QQ$.

\begin{say}\label{say-d=2-order=2}
Let $\tau\in G$ be an element of order $2$.
For the action on $\PP^2$ we have only one possibility:
$\tau: (x_0:x_1:x_2) \longmapsto (-x_0:x_1:x_2)$ and then
$\psi$ has the form $x_0^4+x_0^2 \psi_2(x_1,x_2)+\psi_4(x_1,x_2)=0$,
where $\psi_4$ has no multiple factors (because $B$ is smooth).
For the action on $X$ we have two possibilities:
\begin{eqnarray}
 \tau: &(x_0:x_1:x_2:y) &\longmapsto (-x_0:x_1:x_2:y)\label{eqnarray-d=2-1}
 \\
 \tau: &(x_0:x_1:x_2:y) &\longmapsto (-x_0:x_1:x_2:-y)\label{eqnarray-d=2-2}
\end{eqnarray}
Since $X^\tau$ is an elliptic curve in the case \eqref{eqnarray-d=2-1},
this case does not occur.
Thus we are in the situation of
\eqref{eqnarray-d=2-2}. Then $X^\tau$ consists of four points.
By \eqref{equation-Lefschetz-2} we have $\tr(\tau)=1$. Moreover, $\chi_\tau= \Phi_1^4\Phi_2^3$.
\end{say}
\begin{say}
Assume that $G$ contains an element $\delta$ of order $4$.
Then $\delta^2=\tau$, where $\tau$ is described in \ref{say-d=2-order=2}.
On the other hand, $\chi_\delta=\Phi_4^k\Phi_2^l\Phi_1^m$, where $k>0$.
Then $\chi_\tau=\Phi_2^{2k}\Phi_1^{7-2k}$. 
This contradicts \ref{say-d=2-order=2}.
Thus $G$ does not contain any elements whose order is
divisible by $4$.
\end{say}

\begin{say}\label{say-d=2-order=3}
Let $\theta\in G$ be an element of order $3$.
We have two possibilities for the action on $X$:
\begin{eqnarray}
& \theta: &(x_0:x_1:x_2:y) \longmapsto (\zeta_3 x_0: x_1:x_2:y),\label{eqnarray-d=2-order=3-1}
 \\
&&\psi= x_0^3\psi_1(x_1,x_2)+\psi_4(x_1,x_2),
 \nonumber
 \\[6pt]
 & \theta: &(x_0:x_1:x_2:y) \longmapsto (x_0: \zeta_3 x_1:\zeta_3^2 x_2:y),\label{eqnarray-d=2-order=3-2}
 \\
&&\psi= x_0^4+ a_2x_0^2x_1x_2+ x_0x_1^3+ x_0x_2^3+a_0x_1^2x_2^2.
 \nonumber
\end{eqnarray}
In the case \eqref{eqnarray-d=2-order=3-1} the intersection
$X\cap \{x_0=0\}$ is an elliptic curve of fixed points.
 This contradicts our assumption.

Thus we have case \eqref{eqnarray-d=2-order=3-2}.
Then $X^\theta$ consists of four points and so
$\tr(\theta)=1$. Hence, $\chi_{\theta}=\Phi_1^3\Phi_3^2$.
\end{say}

\begin{say}\label{say-d=2-order=6}
Let $\delta\in G$ be an element of order $6$.
Then $\delta=\tau\theta$, where $\tau$ (resp. $\theta$)
is described in \ref{say-d=2-order=2} (resp. \ref{say-d=2-order=3}).
Hence, $\tr(\delta)=-5$ or $1$.
But in the first case $\Eu(X^\delta)=-2$ and so $\dim X^\delta=1$.
On the other hand, $X^\delta\subset X^\tau$, where $\dim X^\tau=0$.
The contradiction shows that $\tr(\delta)=1$.
\end{say}

\begin{say}\label{say-d=2-order=9}
Let $\delta\in G$ be an element of order $9$.
Since $\chi_\delta$ is divisible by the cyclotomic polynomial
$\Phi_9$, we have $\chi_\delta=\Phi_9\Phi_1$
and so $\tr(\delta)=1$.
The same arguments show that $\tr(\delta)\ge 0$ if $\delta$ is an
element of order $5$ or $7$.
\end{say}

\begin{say}\label{say-d=2-order=15}
Let $\delta\in G$ be an element of order $15$.
As in \ref{say-d=2-order=9} we see that 
$\chi_\delta=\Phi_5\Phi_3\Phi_1$.
Hence $\chi_{\delta^5}=\Phi_3\Phi_1^5$.
This contradicts \ref{say-d=2-order=3}.
\end{say} 
This finishes the proof of Lemma \ref{lemma-d=2}.
\end{proof}

\section{Del Pezzo surfaces of degree $1$}\label{section-d=1}
\begin{say}
Throughout this section, let $X$ be a del Pezzo surface of degree $1$.
Recall that in this case the linear system $|-2K_X|$ determines a double cover $X\to Y\subset \PP^3$,
where $Y$ is a quadratic cone. The corresponding Galois involution $\beta : X\to X$ 
is called the \textit{Bertini involution}.
Its fixed point locus $X^\beta$ is the union of a curve of genus $4$
and a single point $P$.
As in the case $K_X^2=2$, $\beta$ is contained in the center of
$\Aut(X)$ and $-\beta$ acts on $\Pic(X)$ as the reflection with respect to
$Q=K_X^\perp$.

The linear system $|-K_X|$ is an elliptic pencil
whose base locus coincides with $P$ (a single point).
The natural representation $\Aut(X)\to GL(T_{P,X})$ is faithful.
Let $\pi : X \dashrightarrow B=\PP^1$ be the map given by $|-K_X|$.
Here $B$ can be naturally identified with $\PP(T_{P,X})$.
Every singular member $F$ of $|-K_X|$
is an irreducible curve of arithmetic genus $1$.
Hence, $F$ is a rational curve with a unique singularity $R$
which is either a node or a simple cusp.
Computing the topological Euler number we obtain the following.
\end{say}

\begin{stheorem}{\bf Lemma.}\label{lemma-d=1-chi}
Let $\#_{\mathrm {node}}$ \textup(resp. $\#_{\mathrm {cusp}}$\textup)
be the number of nodal \textup(resp. cuspidal rational curves\textup) in the pencil $|-K_X|$.
\[
 \#_{\mathrm {node}}+2\#_{\mathrm {cusp}}=12.
\]
\end{stheorem}

\begin{mtheorem}{\bf Lemma.}\label{lamme-d=1-even}
Any element $\iota\in \Aut(X)$ of order $2$ fixes 
a curve of positive genus.
\end{mtheorem}

\begin{proof}
There are two choices for the action of $\iota$
on $T_{P,X}$: \ $\diag(-1,-1)$ and $\diag(-1,1)$.
In the first case the action coincides with the action
on $T_{P,X}$ of the Bertini involution $\beta$.
Hence, $\iota\comp \beta^{-1}$ acts trivially on $T_{P,X}$
and so $\iota\comp \beta^{-1}$ is the identity map.
In this case, $X^\iota$ contains a curve of genus $4$.
Assume that $\iota$ acts
on $T_{P,X}$ as $\diag(-1,1)$. Then the fixed points locus of $\iota$
contains a smooth curve $C$
passing through $P$ and the action on $B\simeq \PP(T_{P,X})$ is not trivial. Then the restriction
$\pi|_C: C \to B$ cannot be dominant.
Hence $C$ is a fiber of $\pi$ and so $C$ is an elliptic curve.
\end{proof}

\begin{mtheorem}{\bf Lemma.}\label{lemma-d=1-order3}
Let $G=\langle\delta\rangle\subset \Aut(X)$ be a group of order $3$.
Assume that the representation of $G$ in $GL(T_{P,X})$
is given by a scalar matrix.
Then the pair $(X,G)$ is minimal and $X^G$ contains a curve of genus $2$.
\end{mtheorem}

\begin{proof}
Clearly, the action of $\delta$ on $B\simeq\PP(T_{P,X})$ is trivial.
 We claim that 
$X^\delta$ is the union of a smooth irreducible curve $C$ and $P$.
Indeed, if $X^\delta$ contains an isolated point $R\neq P$, then
$\pi$ is well-defined at $R$ and the action of $\delta$ on $T_{R,X}$
in suitable coordinates has the diagonal form $\diag(\zeta_3, \zeta_3^{\pm 1})$.
Let $F=\pi^{-1}(\pi(R))$ be the fiber of $\pi$ passing through $R$.
Since the action on $B$ is trivial,
the differential $d\pi: T_{R,X}\to T_{\pi(R),B}$ is not surjective.
Hence, $R\in F$ is a singular point.
Let $\nu: F'\to F$ be the normalization.
If $R\in F$ is a node, then the cyclic group $G$ has three fixed points $\nu^{-1}(R)$ and $P$
on $F'\simeq \PP^1$, a contradiction.
Hence, $R\in F$ is a cusp. Then
locally near $R$ the map $\nu$ is given by $t \mapsto (t^2,t^3)$.
So the action near $R$ is not free in codimension one.
Again we get a contradiction.

Thus $X^\delta$ consists of $P$ and a smooth curve $C$.
Since $P\not \ni C$, $C$ contains no fibers of $\pi$.
Let $F_1$ be a degenerate fiber of $\pi$. The action
of $G$ on $F_1$ has exactly two fixed
points: $P$ and $R:=\Sing(F_1)$.
Hence, $C\cap F_1=R$ and so $C$ is connected.
Since $C$ is smooth, it must be irreducible.

Denote $r:=\rk \Pic(X)^G$.
By \eqref{equation-Lefschetz-1} and \eqref{equation-Lefschetz-2}
\[
\Eu(X^\delta)= 1+2-2g(C)=3+ \tr(\delta)= 2+ r -\frac 12 (9-r).
\]
The only solution is $r=1$, $g(C)=2$. Then $(X,G)$ is minimal.
\end{proof}

\begin{mtheorem}{\bf Lemma.}\label{lemma-d=1-trivial}
Let $\{1\}\neq G\subset \Aut(X)$
be a group such that the induced action on the pencil $B$ is trivial.
Then some non-identity element of $G$ fixes 
a curve of positive genus.
\end{mtheorem}
\begin{proof}
The group $G$ is contained in the kernel of the composition
\[
G\to GL(T_{P,X})\to PGL(T_{P,X}).
\]
Hence the image of $G$ in $GL(T_{P,X})$ consists of
scalar matrices and so $G$ is a cyclic group.
Let $\delta\in G$ be a generator and let $m>1$ be its order.

The group $G$ acts faithfully on the general member of $|-K_X|$
which is an elliptic curve and $P$ is a fixed point.
Then $G$ must contain an element $\delta$ of order $m=2$ or $3$.
Since the representation $G\to GL(T_{P,X})$ is faithful,
$\delta$ must be either
the Bertini involution $\beta$ or
an element of order $3$ described in Lemma \ref{lemma-d=1-order3}.
The assertion follows.
\end{proof}

\begin{stheorem}{\bf Corollary.}\label{Corollary-d=1-injective}
Let $G\subset \Aut(X)$
be a subgroup such that the natural homomorphism
$G\to \Aut(B)$ is not injective. 
Then some non-identity element of $G$ fixes 
a curve of positive genus.
\end{stheorem}
\begin{proof}
Apply Lemma \ref{lemma-d=1-trivial} to the
kernel of $G\to \Aut(B)$.
\end{proof}

\begin{say}
Now we are ready to finish the proof of Proposition \ref{Proposition-main-del-Pezzo}
in the case $K_X^2=1$.
Assume that any non-identity element of $G$ does not fix a curve of positive genus.
By Corollary \ref{Corollary-d=1-injective} the group $G$ acts faithfully on $B$.
By Lemma \ref{lamme-d=1-even} the order of $G$ is odd.
Hence by the classification of finite subgroups of $PGL_2(\Bbbk)$
(see e.g. \cite{Klein1956}, \cite{Springer1977}) $G$ is a cyclic group.
Let $\delta\in G$ be its generator.
Then the pencil $|-K_X|$ has exactly two invariant
members, say $C_1$ and $C_2$.
We claim that $G$ faithfully acts on $C_1$ and $C_2$.
Indeed, otherwise some non-identity element $\delta\in G$ 
fixes $C_i$ (point-wise). By our assumption $C_i$ has a (unique) singular point,
say $P_i$. Then
$T_{P_i,C_i}=T_{P_i,X}$ and so the action of $G$ on $C_i$ must be faithful,
a contradiction. 
Therefore, $G$ faithfully acts on $C_1$ and $C_2$.

First we assume that both $C_1$ and $C_2$ are smooth
elliptic curves.
Then $G\simeq \ZZ/3\ZZ$ and by Lemma \ref{lemma-d=1-order3}
the element $\delta$ acts on $T_{P,X}$ as $\diag(\zeta_3,\zeta_3^{-1})$.
The fixed points locus
$X^G$ consists of $5$ points $P$, $P_1,\, P_2\in C_1\setminus C_2$
and $P_3,\, P_4\in C_2\setminus C_1$.
Then by \eqref{equation-Lefschetz-2} we have $\tr(\delta)=\tr(\delta^2)=2$
and so $(X,G)$ is not minimal by \eqref{equation-Lefschetz-1}.

Now we assume that $C_1$ has a singular point, say $P_1$.
Since $G$ is cyclic, $P_1$ cannot be an ordinary double point.
Hence, $P_1\in C_1$ is a cusp.
Locally near $P_1$ the normalization is given by $t \mapsto (t^2,t^3)$.
Since the action of $G$ on $X$ is free in codimension one near $P_1$,
the order of $G$ is coprime to $3$.
Then $C_2$ cannot be an elliptic curve, so
 $C_2$ is also a cuspidal rational curve.
 Then $G$
permutes singular members of $|-K_X|$ other than $C_1$ and $C_2$.
By Lemma \ref{lemma-d=1-chi} the order of $G$ divides $12-4=8$, a contradiction.
\end{say}

\section{Conic bundles}\label{section-Conic-bundles}
In this section we consider $G$-surfaces admitting a conic bundle structure.
The convenience of the reader we recall definitions 
and basic facts (see \cite{Dolgachev-Iskovskikh}). 
\begin{say}{\bf Setup.}
Let $X$ be a projective non-singular surface and let $f: X\to B$ be a
dominant morphism, where $B$ is a non-singular curve. We
say that the pair $f$ is a \textit{conic bundle} if
$f_*\OOO_X=\OOO_B$ (i.e. $f$ has connected fibers) and
$-K_X$ is $f$-ample. Then any fiber $X_b$, $b\in B$
is isomorphic to a reduced conic in $\PP^2$.
Let $G$ be a finite group acting on $X$ and $B$.
We say that $f$ is a \textit{$G$-conic bundle} if $f$ is $G$-equivariant.
We say that a $G$-conic bundle $f: X\to B$ is \textit{relatively $G$-minimal} if
$\rk \Pic(X/B)^G=1$.
Throughout this section we assume that $B\simeq \PP^1$ (because $X$ is a rational surface). 
By Noether's formula the number of degenerate fibers equals to $8-K_X^2$.
In particular, $K_X^2\le 8$.

\begin{ppar} 
Moreover, if a $G$-conic bundle $f: X\to \PP^1$ is relatively $G$-minimal,
then $K_X^2\neq 7$. 
From now on $f: X\to B$ denotes a relatively $G$-minimal conic bundle 
with $B\simeq \PP^1$.
If $K_X^2= 8$, then $f$ is a $\PP^1$-bundle, i.e. $X$ is a Hirzebruch surface $\FF_n$.
In this case the action of $G$ on $\Pic(X)$ is trivial and so $H^1(G,\Pic(X))=0$.
For $K_X^2=3$, $5$, and $6$ the pair $(X,G)$ is not minimal:
there exists an equivariant birational morphism to a
$G$-del Pezzo surface $X'$ with $\Pic(X')^G\simeq \ZZ$ and
$K_{X'}^2> K_X^2$ \cite{Iskovskikh-1979s-e}. This case was investigated in the previous 
sections.
\end{ppar}

Thus we have the following
\begin{stheorem}{\bf Proposition.}\label{Proposition-K2=8}
Let $f: X\to \PP^1$ be a $G$-conic bundle with $K_X^2\ge 5$. 
Assume that the surface $X$ is $G$-minimal.
Then $K_X^2=8$ and $X\simeq \FF_n$, where $n\neq 1$. 
Moreover, $X$ is $H^1$-trivial.
\end{stheorem}
\begin{ppar}{\bf Remark.}
Assume that in the  notation of \xref{Proposition-K2=8} the group $G$ is abelian. 
Then it is linearizable if and only if 
it is stably linearizable and if and only if $G$ has a fixed point (see  \cite[\S 8]{Dolgachev-Iskovskikh}
and Lemma \ref{Lemma-Kollar-Szabo})
\end{ppar}
From now on we assume that $K_X^2\le 4$.

\begin{ppar}\label{par-varrho}
Let $G_F$ be the largest group that acts trivially on $B$.
We have an exact sequence
\[
1 \longrightarrow G_F \longrightarrow G \overset{\pi}\longrightarrow G_B \longrightarrow 1,
\]
where $G_B$ acts faithfully on $B$ and $G_F$ acts faithfully on the generic fiber
$X_\eta$. We also have a natural homomorphism
\[
\varrho: G \longrightarrow \Aut(\Pic(X)).
\]
Since $B\simeq \PP^1$ and $K_X^2\le 5$, the group $\ker(\varrho)$
fixes point-wise any section with negative self-intersection.
In particular, this implies that $\ker(\varrho)\subset G_F$ and $\ker(\varrho)$
is a cyclic group.
\end{ppar}

\begin{ppar} {\bf Notation.}\label{conic-bundle-notation}
 Let $f: X\to B\simeq \PP^1$ be a relatively $G$-minimal $G$-conic bundle and let $F$ be a typical fiber. Let
$F_1,\dots, F_m$ be all the degenerate fibers, let $R_i$ be the singular point of $F_i$, and let 
$P_i:=f(F_i)$.
Thus, $F_i=f^{-1}(P_i)=F_i'+F_i''$ and $F_i'\cap F_i''=\{R_i\}$.
Let $\Delta:=\{P_1,\dots, P_m\}$ be the discriminant locus.
\end{ppar}
\end{say}

\begin{mtheorem}{\bf Lemma (cf. \cite[Lemmas 3.9-3.10]{Blanc2011}).}
\label{Corollary-cyclic}
In the notation of \xref{conic-bundle-notation}
assume that any non-identity element of $G$ 
does not fix a curve of positive genus.
Let $\delta\in G$ be an element of order $n>1$.
Then one of the following holds:
\begin{enumerate}
\item
$\delta$ does not switch components of any degenerate fiber, 
\item
there are exactly two degenerate fibers whose components are switched by $\delta$, or
\item
$\delta$ switches components of exactly one degenerate fiber, 
say $F_1$. In this case, 
$\delta^{2}$ acts on $B$ trivially and $\delta$ acts on $B$ non-trivially.
Moreover, $\delta^{2}$ 
switches components of exactly two degenerate fibers \textup(other than $F_1$\textup).
\end{enumerate}
\end{mtheorem}
\begin{proof}
Let $F_1,\dots, F_r$ be all the degenerate fibers whose
components are switched by $\delta$. We assume that $r>0$
(otherwise we are in the situation of (i)).

First we consider the case where the action of $\delta$ on $B$ is 
trivial. Then $\delta$ has exactly two fixed points on any smooth fiber.
Hence, $X^\delta$ contains a (smooth) curve $C$. 
For $i\in \{1,\dots, r\}$, each intersection point $C\cap F_i$ is a single point 
which must coincide with 
$R_i=\Sing(F_i)$.
So, $C$ is connected and the ramification locus of 
the double cover $f_C: C\to B$ coincides with $\{P_1,\dots, P_r\}$.
In particular, $r$ is even.
If $r>2$, then $C$ is a curve of genus $(r-2)/2>0$, a contradiction.
Hence, $r=2$. 

Now consider the case where the action of $\delta$ on $B$ is 
non-trivial.
Since $\delta$ has exactly two fixed points on $B$, we have $r\le 2$.
Assume that $r=1$. 
If any element of the group $\langle\delta \rangle$ does not switch
components of any fiber except for $F_1$, then we can run a relative
$\langle \delta\rangle$-minimal model program on $X$ so that 
the resulting surface has a relatively 
$\langle\delta \rangle$-minimal conic bundle structure over $B$ 
with exactly one degenerate fiber. 
It is easy to see (see e.g. \cite[Lemma 5.1]{Dolgachev-Iskovskikh}) that 
this is impossible.
Hence some element $\delta^k$, where $k>1$, switches
components of a fiber $F_2\neq F_1$.
Take $k$ to be minimal possible.
The points $f(F_2)$ and $f(F_1)$ are fixed by $\delta^k$.
By our assumption $r=1$, the point $f(F_2)$ is not fixed by $\delta$.
This is possible only if $\delta^k$ acts trivially on $B$.
According to the above considered case, $\delta^k$ switches 
components of exactly two fibers, so the $\langle\delta \rangle$-orbit of 
$F_2$ consists of two elements. Therefore, $k=2$. 
\end{proof}

Now we are going to classify  $H^1$-trivial  $G$-conic bundles
with $K_X^2\le 4$.
There are two essentially different cases: $\ker(\varrho)=\{1\}$ and $\neq \{1\}$.

\subsection*{Case $\ker(\varrho)=\{1\}$.}

\begin{mtheorem}{\bf Theorem.}\label{theorem-conic-bundle-1}
Let $f: X\to B=\PP^1$ be a relatively $G$-minimal $G$-conic bundle
with $K_X^2\le 4$.
Assume that
$(X,G)$ is $H^1$-trivial and $\ker(\varrho)=\{1\}$. Then 
$G\simeq \tilde\DD_n$, where $n=6-K_X^2$ is odd,
$G_F\simeq \ZZ/2\ZZ$ is the center of $G$, $G/G_F\simeq \DD_n$,
and the action is given by Construction \xref{Construction-conic-bundle-non-exceptional}
below\footnote{For $n=5$ see also \cite[Th. 6.5]{Tsygankov2011e}}.
\end{mtheorem} 

\begin{ppar} {\bf Remark.}
In the case $n=3$ the surface $X$ is not $G$-minimal: 
contracting an invariant horizontal $(-1)$-curve we get a
quartic del Pezzo surface (see 
\eqref{equation-d=4}-\eqref{equation-d=4-action}
and Remark \ref{Remark-d=4-case-realized}).
\end{ppar}

\begin{say} {\bf Construction (cf. \cite[5.12]{Dolgachev-Iskovskikh}, \cite[3.2]{Tsygankov2011e}).}
\label{Construction-conic-bundle-non-exceptional}
Let $n\ge 3$ be an odd integer.
The representation \eqref{equation-binary-dihedral-group}
induces a faithful action 
$\varsigma_1: \DD_{n} \longrightarrow \Aut(\PP^1)$.
Consider another faithful action $\varsigma_2: \DD_{n} \longrightarrow \Aut(\PP^1)$
\[
\tilde r \mapsto
\begin{pmatrix}
\zeta_{n}& 0
\\
0& \zeta_{n}^{-1}
\end{pmatrix}
\qquad
\tilde s \mapsto
\begin{pmatrix}
0_{\phantom{n}}& -1_{\phantom{n}}
\\
-1_{\phantom{n}}& 0_{\phantom{n}}
\end{pmatrix}
\]
Clearly we have $\lambda\comp \varsigma_1=\varsigma_2 \comp \lambda$,
where the map $\lambda: \PP^1\to \PP^1$ is given by $\lambda:x \mapsto x^2$.
Consider also the action 
\[
 \varsigma=\varsigma_1\times \varsigma_2: \DD_{n} \longrightarrow \Aut(\PP^1\times \PP^1).
\]
The curves 
\begin{eqnarray*}
\Gamma&:=&\{ (x,y) \in \PP^1\times \PP^1 \mid x^2=y\}, 
\\
L&:=&\{(x,y) \in \PP^1\times \PP^1 \mid y^n=1\} 
\end{eqnarray*}
are $\DD_{n}$-invariant. Let $L_k:= \{(x,y) \mid y=\zeta_n^k\}$
be a component of $L$. It is easy to see that $L_k$ meets $\Gamma$ transversally
at two points. Now we explicitly construct a double cover $\pi:Y\to \PP^1\times \PP^1$ 
branched over $\Gamma+L$. In homogeneous coordinates on $\PP^1\times \PP^1$ 
the curve $\Gamma+L$ is given by
\[
\phi:=(x_1^2y_0-x_0^2y_1)(y_1^n-y_0^n)=0.
\]
For short, we put $q:=(n+1)/2$. 
Let $\nu:\PP^1\times \PP^1 \longrightarrow \PP^{n+2}$ be the Segre embedding
\begin{gather*}
\nu:\left((x_0:x_1), (y_0,y_1)\right) \longmapsto \left(t_{0,0}, \dots, t_{0,q},t_{1,0}, \dots, t_{1,q}\right),
\quad \text{where}\\
t_{a,b}= x_0^{1-a}x_1^{a} y_0^{q-b}y_1^{b},\qquad 0\le a\le 1,\ 0\le b\le q.
\end{gather*}
Clearly, $\phi$ can be written as a homogeneous polynomial of degree $2$ in the $t_{a,b}$'s.
Thus we can exhibit $Y\subset \PP^{n+3}$ as the intersection of the hypersurface 
\[
z^2=\phi(t_{0,0},\dots,t_{1,q})
\]
with the projective cone which is the preimage of $\nu(\PP^1\times \PP^1)$ under the projection 
\[
\PP^{n+3}\dashrightarrow \PP^{n+2}\supset \nu(\PP^1\times \PP^1),
\quad 
\left(z,t_{0,0},t_{0,1},\dots\right) \longmapsto \left(t_{0,0},t_{0,1},\dots \right).
\]
Let $\sigma :\DD_n\to \{\pm 1\}$ be as in \ref{notation-1}.
Consider the group 
\[
\{ (\delta,\alpha) \in \DD_n \times \langle\zeta_4\rangle \mid \sigma(\delta)=\alpha^2\}. 
\]
This group is a non-trivial central extension of $\DD_n$ by 
$\ZZ/2\ZZ$ and it is isomorphic to $\tilde \DD_n$. 
By the above construction 
 we see that $\tilde \DD_n$ acts on $Y$ so that $\pi$ is equivariant. 
The projection of $\PP^1\times \PP^1$ to the second factor induces a rational curve fibration
$Y\to \PP^1$ whose fibers are irreducible except for those corresponding to two ramification points
of the double cover $\Gamma \to \PP^1$. 
Let $\bar L_k:=\pi^{-1}(L_k)$. 
There are exactly $2n$ nodes 
$Q_1',Q_1'',\dots,Q_n',Q_n''\in Y$, where $\{Q_k',Q_k''\}=\pi^{-1}(\Gamma\cap L_k)$.
Let $\tilde Y\to Y$ be the minimal resolution and let $\tilde Y\to X$
the contraction of the $\tilde L_k$'s, the proper transforms of the $\bar L_k$'s.
Then $f: X \to \PP^1$ is a $\tilde \DD_n$-conic bundle with $n+2$ degenerate fibers fitting to the following
commutative diagram:
\begin{equation}
\begin{array}{c}
\label{equation-cd-1}
\xymatrix@R17pt{
\tilde Y\ar[d]\ar[rr]&&Y\ar[d]^{\pi}
\\
X&&\PP^1\times \PP^1\ar[ll]
} \end{array}
\end{equation}
\end{say}

\begin{proof}[Proof of Theorem \xref{theorem-conic-bundle-1}]
Assume that $\varrho$ is injective.
Then so $\varrho_F: G_F\to \Aut(\Pic(X))$ is.

\begin{stheorem}{\bf Lemma.}
 $G_F\neq\{1\}$.
\end{stheorem}

\begin{proof}
Indeed, otherwise $G$ faithfully acts on $B=\PP^1$.
For any degenerate fiber $F_i$, there exist an element 
$\delta \in G$ switching the components of $F_i$.
In particular, $\ord(\delta)=2k$ for some $k$.
Clearly, we may assume that $k=2^l$.
By Lemma \ref{Corollary-cyclic} there exists exactly one 
more degenerate fiber $F_j\neq F_i$ whose components
are switched by $\delta$.
Thus $X^\delta=\{R_i,\, R_j\}$. If $k=1$, then the holomorphic Lefschetz fixed point formula
implies that the cardinality of $X^\delta$ equals to $4$, a contradiction.
Hence, $k>1$. Put $\gamma:=\delta^{k}$. 
It is easy to see that 
$X^\gamma= F_i^\gamma\cup F_j^\gamma$.
Since $X^\gamma$ is $\delta$-invariant and smooth,
one can see that it is zero-dimensional and 
consists of exactly $6$ points. Again we get a contradiction
by the holomorphic Lefschetz fixed point formula.
This proves our lemma.
\end{proof}

The group $G_F$ interchanges pair-wise components of (some) degenerate fibers.
So, there exists an embedding 
\[
G_F \hookrightarrow \Sym_2\times
\cdots \times \Sym_2.
\]
On the other hand, $G_F$ acts faithfully on a typical fiber,
so there exists an embedding $G_F \hookrightarrow PGL_2(\Bbbk)$.
This immediately gives us either $G_F\simeq \ZZ/2\ZZ$
or $G_F\simeq (\ZZ/2\ZZ)^2$ 
(see \cite[Th. 5.7]{Dolgachev-Iskovskikh}).

Consider the case $G_F\simeq (\ZZ/2\ZZ)^2$.
Then $G_F=\{1,\, \tau_1,\, \tau_2,\, \tau_3\}$,
where the $\tau_j$'s are distinct elements of order $2$.
Fix $i\in \{1,\dots,m\}$.
The point $R_i$ is fixed by $G_F$.
The actions of all the $\tau_j$'s on $T_{R_i,X}$ cannot have the (same) form $\diag(-1,-1)$.
Hence at least one of them, say $\tau_1$, is of type $\diag(1,-1)$
(in suitable coordinates).
Then $\tau_1$ must switch the components of $F_i$.
Indeed, otherwise $\tau_1$ fixes point-wise a component of $F_i$.
But this is impossible because $\tau_1$ acts trivially on $B$.
Moreover, for each singular fiber $F_i$, exactly two elements of $G_F$ 
switch the components of $F_i$.
Taking Lemma \ref{Corollary-cyclic} into account, we see that 
$\Delta$ consists of three elements.
This contradicts our assumption $K_X^2\le 4$.

Therefore, $G_F\simeq \ZZ/2\ZZ$.
Let $\tau \in G_F$ be the element of order $2$.
Since $\varrho(\tau)\ne \id$, 
by Lemma \ref{Corollary-cyclic} the element $\tau$ switches
components of exactly two degenerate fibers, say $F_{r-1}$ and $F_r$.
By our assumption $K_X^2\le 4$, we have $r>2$.
Then the set $\{P_{r-1},\, P_r\}$ is $G_B$-invariant.
This is possible only if $G_B$ is either cyclic or dihedral.
Let $C$ be the one-dimensional part of $X^\tau$.
As in the proof of Lemma \ref{Corollary-cyclic} we see that
$C$ is a smooth rational curve and $f_C: C\to B$ is a double cover ramified 
over $\{P_{r-1},\, P_r\}$. The group $G_B=G/G_F$ faithfully acts on $C$
so that $f_C$ is $G_B$-equivariant.

Let $\delta\in G$ be an element that switches
the components of $F_1$. If $\delta$ does not permute $F_{r-1}$ and $F_r$, then
$\delta$ fixes three points $P_{r-1},\, P_r,\, P_1\in B=\PP^1$.
So, it trivially acts on $B$, that is, $\delta\in G_F$, a contradiction.
Thus $\delta$ permutes $F_{r-1}$ and $F_r$.
Let $\upsilon\in \Aut(C)$ be the Galois involution of $f_C$
and let $G_C\subset \Aut(C)$ be the (isomorphic) image of $G_B$.
Since $G_B$ faithfully acts on $B$, $\upsilon\notin G_C$.
On the other hand, $\upsilon$ commutes with any element of $G_C$.
Hence, $G_C$ and $\upsilon$ generate a subgroup $G_C'=G_C\times \langle \upsilon\rangle\subset \Aut(C)$
so that the set $\{R_{r-1},\,R_r\}\subset C$
is $G_C'$-invariant. By the classification of finite subgroups of
$\Aut(\PP^1)$ we see that $G_C'\simeq \DD_{2n}$, where 
$n$ must be odd (because $\upsilon\notin \DD_{n}\subset \DD_{2n}$).
In particular, $G_B\simeq \DD_n$.
For $i=1,\dots,r-2$ we have $C\cap F_i'=\{R_i'\}$ and $C\cap F_i''=\{R_i''\}$,
where the points $R_i'$ and $R_i''$ are permuted by $\upsilon$
and have non-trivial stabilizers in $G_C$.
There are only three non-trivial orbits of $\DD_{2n}$ on $C\simeq\PP^1$:
$O_{2n}$, $O_{2n}'$ and $O_2$ (\cite{Klein1956}, \cite{Springer1977}).
They have $2n$, $2n$, and $2$ elements, respectively.
Since $\upsilon$ can not fix any element of $O_{2n}$ and $O_{2n}'$,
we may assume that $O_{2n}'$ form one $\DD_{n}$-orbit
and $O_{2n}$ splits in the union of two $\DD_{n}$-orbits.
Then $O_{2n}$ coincides with $C\cap (\cup_{i=1}^{r-2} F_i)$
and so $n=r-2$. 
Recall that $n$ is odd and $G$ is a central extension of $G_B\simeq \DD_n$ by $G_F\simeq \ZZ/2\ZZ$.
We claim that $G\simeq \tilde \DD_{n}$.
Indeed, otherwise 
$G=G_B\times G_F\simeq \DD_n\times \ZZ/2\ZZ$.
Take $\delta$ as above. Then $\delta$ fixes $P_1$.
Since $G\simeq \DD_n\times \ZZ/2\ZZ$, we have 
$\ord(\delta)=2$.
The action of $\delta$ on $T_{R_1,X}$ has the form 
$\diag(1,-1)$. Hence, $\delta$ fixes point-wise a 
(smooth) curve $D$ passing through $R_1$. 
Since $\delta$ switches the components of $F_1$, $D$ is not a component of $F_1$.
Hence, $D$ dominates 
$B$ and $\delta\in G_F$, a contradiction. Thus $G\to G_B$ is not split and 
so $G\simeq \tilde \DD_n$.

Now we construct the following $G$-equivariant commutative diagram
\begin{equation}
\label{equation-cd-2}
\begin{array}{c}
\xymatrix{
\tilde Y\ar[d]\ar@/^12pt/[rrr]&&Z\ar[dr]^{\upsilon}\ar[dl]_{\mu}&Y\ar[d]^\pi
\\
X\ar[r]&X/\langle\tau\rangle&&\FF_e
} 
\end{array}
\end{equation}
Here $X/\langle\tau\rangle$ has $n=r-2$ nodes which are images of $R_1,\dots,R_{n}$,
$\mu$ is the minimal resolution and $\upsilon$ is the contraction of the proper transforms 
of $R_1',R_1'',\dots, R_{n}',R_{n}''$. It is easy to see that the image of $\upsilon$ must be a smooth
geometrically ruled surface.
On the other hand, to arrive to $\FF_e$ from $X$ we can blowup the points $R_1,\dots, R_{n}$ first.
We get $\tilde Y$. The action of $G$ lifts to $\tilde Y$ and $\tilde Y\to Y\to \FF_e$ 
is the Stein factorization. 
Let $E_1,\dots,E_{n}$ be $\mu$-exceptional divisors and let $L_k:=\upsilon(E_k)$.
Let $C_\bullet\subset \FF_e$ be the proper transform of 
$C/\langle\tau\rangle\subset X/\langle\tau\rangle$.
Clearly, $\pi$ is a double cover branched over $C_\bullet+L_1+\cdots+L_{n}$. 
Comparing \eqref{equation-cd-2} and \eqref{equation-cd-1} we see that
it remains to show that $e=0$, i.e. $\FF_e\simeq \PP^1\times \PP^1$. 
We can write 
$C_\bullet\sim 2s+a F_\bullet$, where $s$ is the minimal section and 
$F_\bullet$ is a fiber of $\FF_e$. Since $C_\bullet$ is an irreducible 
smooth rational curve, we get two possibilities:
$(e,a)=(0,1)$ and $(1,2)$.
Since the branch divisor $C_\bullet+L_1+\cdots+L_{n}$ is divisible 
by $2$ and $n$ is odd, we see that the second case is impossible.
This proves Theorem \xref{theorem-conic-bundle-1}.
\end{proof}

\subsection*{Case $\ker(\varrho)\neq \{1\}$.}

\begin{say} {\bf Definition \cite{Dolgachev-Iskovskikh}.}\label{Definition-exceptional}
A conic bundle $f: X\to \PP^1$ is said to be \textit{exceptional} 
if for some positive integer $g$ the number of degenerate fibers equals to
$2g+2$ and 
there are two disjoint sections $C_1$ and $C_2$ with $C_1^2=C_2^2=-(g+1)$.
\end{say}

\begin{mtheorem}{\bf Theorem.}\label{theorem-conic-bundle-2}
Let $f: X\to \PP^1$ be a relatively $G$-minimal $G$-conic bundle
with $K_X^2=6-2g\le 4$.
Assume that
$(X,G)$ is $H^1$-trivial and $\ker(\varrho)\neq \{1\}$.
Then we have

\begin{enumerate}
 \item \label{theorem-conic-bundle-0}
$f$ is exceptional, in particular, $K_X^2$ is even;
\item \label{theorem-conic-bundle-i}
$G_F=\ker(\varrho)$ and it is a non-trivial cyclic group;

\item \label{theorem-conic-bundle-b}
either $G_B\simeq \DD_{n}$ or $G_B\simeq \Sym_4$;
\item \label{theorem-conic-bundle-iv}
the action of $G$ on $X$ is given by the Construction \xref{Construction-conic-bundles-exceptional}
below.
\end{enumerate}
\end{mtheorem}

 The following is 
a particular case of the general construction 
\cite[\S 5.2]{Dolgachev-Iskovskikh}.

\begin{say} {\bf Construction \cite[\S\S 5.2-5.3]{Dolgachev-Iskovskikh}.}
\label{Construction-conic-bundles-exceptional}
First we fix some data.
Let $\tilde G_B\subset SL_2(\Bbbk)$ be a finite non-cyclic subgroup
and let $G_B=\tilde G_B/\{\pm \id\}$ be its image in $PSL_2(\Bbbk)$.
Fix two homomorphisms $\varsigma, \chi_B: G_B\to \{\pm 1\}$, where $\chi_B$ is surjective
(we assume that such a homomorphism $\chi_B$ exists).
We also regard $\varsigma$ and $\chi_B$ as characters defined on $\tilde G_B$.
Let $g\ge 1$ and let $Y$ be
the hypersurface in $\PP(g+1,g+1,1,1)$ given by $x_1x_2=\psi(y_1,y_2)$,
where $\psi(y_1,y_2)$ is a homogeneous $\tilde G_B$-semi-invariant 
of degree $2g+2$ and weight $\varsigma$. Thus 
$\delta(\psi)= \varsigma(\delta)\psi$ for all $\delta\in \tilde G_B$.
We assume also that $\psi$ has no multiple factors.
Put 
\[
\Gamma:=\{ (h,\delta)\in GL_2(\Bbbk)\times \tilde G_B\ \mid \
h(x_1x_2)= \varsigma(\delta) x_1x_2\}.
\]
It is easy to see that $\Gamma$ naturally acts on $Y$
and the kernel of the action coincides with
\[
K:= 
\left\langle \left((-1)^{g+1}\id,-\id\right) \right\rangle.
\]
Thus $\Aut(Y)\supset \Gamma/K$.
Denote by $p:\Aut(Y) \to G_B$ the homomorphism induced by the projection to the second factor. 
The surface $Y$ has two singular points which are of type $\frac{1}{g+1}(1,1)$.
Let $X\to Y$ be the minimal resolution. The projection 
$(x_1:x_2:y_1:y_2) \dashrightarrow (y_1:y_2)$ induces a conic bundle structure
$f: X \to \PP^1=B$ whose degenerate fibers correspond to the zeros of
$\psi$. In particular, $K_X^2=6-2g$.

The action on the set $\Sing(Y)=\{(1:0:0:0),\, (0:1:0:0)\}$
defines a homomorphism $\chi: \Aut(Y)\to \{\pm 1\}$.
Now, take a subgroup $G\subset \Gamma/K$ such that the restriction $\chi_G: G\to \{\pm 1\}$
and the projection $p_G: G\to G_B$ are surjective, and $\ker(p)\cap G\subset \ker(\chi)$.
Thus $\chi$ descends to a character $\chi_B: G_B\to \{\pm 1\}$.

There are the following possibilities:
\renewcommand{\arraystretch}{1.5}
\setlength{\tabcolsep}{15pt}
\newcounter{nom}
\renewcommand{\thenom}{\rm\arabic{nom}$^o$}
\def\bnom{\begingroup\refstepcounter{nom}\rm\arabic{nom}$^o$}
\def\enom{\endgroup}

\par\medskip\noindent
\begin{tabular}{p{3pt}|p{18.7pt}|p{17pt}|l|l|l}
{\rm No.} & $g$ & $G_B$ & $\psi$&$\varsigma$&$\chi_B$
 \\\hline
\bnom\label{theorem-conic-bundle-1-S4-1}\enom&$2$ &$\Sym_4$&$\psi_6$& $\sgn$ & $\sgn$ 
 \\
\bnom\label{theorem-conic-bundle-1-S4-2}\enom&$5$& $\Sym_4$&$\psi_{12}$&$1$& $\sgn$ 
 \\
\bnom\label{theorem-conic-bundle-1-S4-3}\enom& $8$& $\Sym_4$&$\psi_{6}\psi_{12}$&$\sgn$& $\sgn$ 
 \\
\bnom\label{theorem-conic-bundle-1-Dn-i}\enom& $\ge 1 $ &$\DD_{2g+2}$& $y_1^{2g+2}-y_2^{2g+2}$&$\xi\cdot \sigma^{g-1}$&$\sigma$ or $\xi$
 \\
\bnom\label{theorem-conic-bundle-1-Dn-ii}\enom& $\ge 1 $ &$\DD_{g+1}$&$y_1^{2g+2}-y_2^{2g+2}$&$\sigma^{g}$&$\sigma$
 \\
\bnom\label{theorem-conic-bundle-1-Dn-iii}\enom& $\ge 1 $ &$\DD_{2g}$&$y_1y_2(y_1^{2g}-y_2^{2g})$&$\xi\cdot \sigma^{g-1}$&$\xi$
\end{tabular}
\par\medskip\noindent
where $\psi_6=y_1y_2(y_1^{4}-y_2^{4})$ and $\psi_{12}=y_1^{12}-33y_1^8y_2^4-33y_1^4y_2^8+y_2^{12}$,
and, for even $n$, the homomorphism $\xi: \DD_n\to \{\pm 1\}$ is defined by $\xi(r)=-1$, $\xi(s)=-1$.
\end{say}

\begin{proof}[Proof of Theorem \xref{theorem-conic-bundle-2}\xref{theorem-conic-bundle-0}]
Since $\ker(\varrho)\neq \{0\}$, the conic bundle $f$ is exceptional by 
\cite[Proposition 5.5]{Dolgachev-Iskovskikh}. 
In particular, we can write $m=2g+2$, where $g\in \ZZ_{>0}$.
\end{proof}
Let $C_1$ and $C_2$ are disjoint $-(g+1)$-sections (see Definition \ref{Definition-exceptional}).
\begin{proof}[Proof of Theorem \xref{theorem-conic-bundle-2}\xref{theorem-conic-bundle-i}]
Recall that $\ker(\varrho)\subset G_F$ (see \ref{par-varrho}).
If there exists an element $\delta\in G_F$ that switches
$C_1$ and $C_2$, then $\delta$ switches components of all degenerate fibers.
Since the number of degenerate fibers equals to
$2g+2\ge 4$, this contradicts Lemma \ref{Corollary-cyclic}.
Hence both $C_1$ and $C_2$ are $G_F$-invariant
and then any component of a degenerate fiber also must be $G_F$-invariant.
Since $K_X$ and the components of the fibers generate 
a subgroup of index $2$ in $\Pic(X)$,
we have  $G_F=\ker(\varrho)$.
Finally, the action of $G_F$ on a typical fiber $F$
has two fixed points $C_1\cap F$ and $C_2\cap F$.
Then $G_F$ must be cyclic.
\end{proof}

\begin{stheorem}{\bf Corollary.} \label{Corollary-chi-homomorphism}
Let $\chi: G\to \{\pm 1\}$ be the \textup(surjective\textup) homomorphism
induced by the action on $\{C_1,\, C_2\}$. Then $G_F\subset \ker(\chi)$.
Thus $\chi$ passes through a surjective homomorphism $\chi_B: G_B\to \{\pm 1\}$.
\end{stheorem}

\begin{proof}[Proof of Theorem \xref{theorem-conic-bundle-2}\xref{theorem-conic-bundle-b}]
Suppose that $G_B$ is cyclic.
By \xref{theorem-conic-bundle-i} of our theorem $G_B\neq \{1\}$.
Thus $G_B$ has exactly two fixed points $P',\, P''\in B$
and acts freely on $B\setminus \{P',\, P''\}$.
For any degenerate fiber $F_i$ there exists an element
$\delta\in G$ that switches components of $F_i$.
Then $P_i=f(F_i)$ must coincide with $P'$ or $P''$.
Hence, $f$ has at most two degenerate fibers, a contradiction.
Thus $G_B$ is not cyclic.

Recall that $G_B\subset PGL_2(\Bbbk)$.
By the classification of finite subgroups in $PGL_2(\Bbbk)$
(see e.g. \cite{Klein1956}, \cite{Springer1977})
we have $G_B\simeq \DD_{n}$,
$\mathfrak A_4$, $\Sym_4$, or $\mathfrak A_5$.
By Corollary \ref{Corollary-chi-homomorphism} we have
$G_B\not\simeq \mathfrak A_4$, $\mathfrak A_5$.
\end{proof}

\begin{stheorem}{\bf Lemma.} \label{Lemma-in}
For $P_i=f(F_i)$, let $G_i\subset G_B$ be its stabilizer.
Then $G_i$ is a cyclic group generated by an element $\tau_i$ such that
$\chi_B(\tau_i)=-1$.
\end{stheorem}
\begin{proof}
Since the representation of $G_i$ on $T_{P_i,B}$ is faithful, $G_i$ is cyclic.
The components of $F_i$ are switched by some element $\delta_i\in G$.
Then $\chi(\delta_i)=-1$ and the image of $\delta_i$ is contained in $G_i$.
\end{proof}
\par\medskip\noindent
\textit{Proof of Theorem \xref{theorem-conic-bundle-2}\xref{theorem-conic-bundle-iv}.}
Basically, this is the third construction of exceptional conic bundles
in \cite[5.2]{Dolgachev-Iskovskikh}. We have to prove only
\ref{theorem-conic-bundle-1-S4-1}-\ref{theorem-conic-bundle-1-Dn-iii}.

\begin{say}
Define a homogeneous semi-invariant $\psi(y_1,y_2)$ so that it
vanishes at $P_1$,\dots, $P_{2g+2}\in \PP^1_{y_1,y_2}$ with multiplicity one
and does not vanish everywhere else.

\begin{stheorem}{\bf Lemma.}
Let $G_i\subset G_B$ be the stabilizer of 
 $P_i=f(F_i)$ and let $\tau_i$ be its generator.
Then the set $\Delta:=\{P_1,\dots, P_{2g+2}\}$
satisfies the following property:
\begin{itemize}
\item
the fixed point locus 
 $B^{\tau_i}$ is contained in $\Delta$.
\end{itemize}
In particular, $\Delta$ is the union of some non-trivial $G_B$-orbits.
\end{stheorem}
\begin{proof}

Let $\hat\tau_i\in G$ be a preimage of $\tau_i$. 
By construction, $\hat\tau_i$ switches components of $F_i$.
If $F_i$ is the only fiber whose components are switched by 
$\hat\tau_i$, then $\hat\tau_i$ is as in Lemma \ref{Corollary-cyclic}(iii).
But then 
$\hat\tau_i^2\in G_F=\ker(\varsigma)$ and so $\hat\tau_i^2$ 
does not switch components of any fiber.
This contradicts Lemma \ref{Corollary-cyclic}(iii).
Hence $\hat\tau_i$ switches the components of two fibers: $F_i$ and $F_j\neq F_i$.
Therefore, $B^{\tau_i}=\{f(F_i),\, f(F_j)\}\subset \Delta$
\end{proof}
\end{say}

\begin{say}
Consider the case $G_B\simeq \Sym_4$.
Then $\chi$ coincides with the sign map $\sgn: \Sym_4\to \{\pm 1\}$.
There are only three non-trivial orbits of $\Sym_4$ on $\PP^1$:
$O_{12}$, $O_{8}$, and $O_{6}$
(see e.g. \cite{Klein1956}, \cite{Springer1977}).
They have $12$, $8$, and $6$ elements, respectively.
The corresponding semi-invariants have the form
$\psi_{12}=y_1^{12}-33y_1^8y_2^4-33y_1^4y_2^8+y_2^{12}$,\
$\psi_8=y_1^{8}+14y_1^4y_2^4+y_2^{8}$, and
$\psi_6=y_1y_2(y_1^{4}-y_2^{4})$.
By Lemma \ref{Lemma-in} for any point $P_i\in \Delta$
its stabilizer $G_i\subset G_B$ is generated by an odd permutation.
So, the order of $G_i$ equals to $2$ or $4$ and $O_8\not\subset \Delta$.
Hence there are the following possibilities: $\Delta=O_{12}$,
$\Delta=O_{6}$, and $\Delta=O_{6}\cup O_{12}$.
\end{say}

\begin{say}
Now consider the case $G_B\simeq \DD_n$.
We use the presentation \eqref{equation-dihedral-group}.
There are only three non-trivial orbits of $\DD_n$ on $\PP^1$:
$O_n$, $O_n'$ and $O_2$ (\cite{Klein1956}, \cite{Springer1977}).
They have $n$, $n$, and $2$ elements, respectively.
The corresponding semi-invariants of $\DD_n$
have the form $\psi_n=y_1^n-y_2^n$, $\psi_n'=y_1^n+y_2^n$, $\psi_2=y_1y_2$.
Since $\Delta$ contains at least four points, $\Delta\neq O_2$.
Thus we may assume that $O_n\subset \Delta$.
Assume that $\Delta\supset O_n\cup O_n'$. Then
any element $\tau\in \DD_n\setminus \langle r\rangle$
generates the stabilizer of some point $P_i\in\Delta$.
By Lemma \ref{Lemma-in} the character $\chi$ takes value $-1$ on $\DD_n\setminus \langle r\rangle$.
Hence, $\chi(r)=1$,\ $r$ cannot generate the stabilizer of a point of $\Delta$
and so $O_2\not\subset \Delta$.
Thus for $\Delta$ we have the following possibilities:
$\Delta= O_n$, $O_n\cup O_n'$, and $O_n\cup O_2$,
corresponding to
\ref{theorem-conic-bundle-1-Dn-i},
\ref{theorem-conic-bundle-1-Dn-ii},
and \ref{theorem-conic-bundle-1-Dn-iii},
respectively.
Finally, $\chi_B$ can be computed by using Lemma \ref{Lemma-in}.
This proves Theorem \xref{theorem-conic-bundle-2}. \hspace{\fill} $\square$
\end{say}

\begin{mtheorem}{\bf Corollary.}\label{Corollary-Iskovskikh-surface}
Let $f: X\to B=\PP^1$ be a relatively $G$-minimal $G$-conic bundle, where $G$ is an abelian group.
Assume that $f$ has at least one degenerate fiber,
$(X,G)$ is $G$-minimal and $H^1$-trivial. Then the following assertions hold:
\begin{itemize}
\item
$K_X^2=4$, $G\simeq \ZZ/4\ZZ\oplus \ZZ/2\ZZ$, $f$ has exactly $4$ degenerate fibers,
\item
the image of $G$ in $\Aut(B)$ is isomorphic to $\ZZ/2\ZZ\oplus \ZZ/2\ZZ$, and
$f$ is an exceptional conic bundle with $g=1$.
\item
There are two disjointed sections $C_1$ and $C_2$ which are $(-2)$-curves.
Moreover, $X$ is a weak del Pezzo surface, that is, $-K_X$ is nef and big.
\item
The anti-canonical model $\bar X\subset \PP^4$ an intersection of two quadrics
whose singular locus consists of two ordinary double points and the line joining
them does not lie on $\bar X$.
\end{itemize}
\end{mtheorem}

\begin{say}{\bf Remark.}\label{Remark1-Iskovskikh-surface}
The surface $X$ and group $G$ described above are extremal in many senses.
According to \cite[\S 7]{Blanc2009} \ $G$ is the only finite
abelian subgroup of $\Cr_2(\Bbbk)$ which is not conjugate to a group of
automorphisms of $\PP^2$ or $\PP^1\times \PP^1$ but whose non-trivial elements do not fix any
curve of positive genus.
The intersection of two quadrics $\bar X\subset \PP^4$ as above is called
the \textit{Iskovskikh surface} \cite{Coray-Tsfasman-1988}.
This is the only intersection of two quadrics in $\PP^4$
for which the Clean Hasse Principle can fail
\cite{Iskovskikh1971-e},
\cite{Coray-Tsfasman-1988}.
\end{say}
\begin{say}{\bf Remark.}\label{Remark2-Iskovskikh-surface}
In the notation of Corollary \ref{Corollary-Iskovskikh-surface},
it is easy to see that the group $G=\ZZ/4\ZZ\oplus \ZZ/2\ZZ$ has no any fixed points on 
$X$. 
Hence $(X,G)$ is not stably linearizable (see Lemma \ref{Lemma-Kollar-Szabo}). 
Moreover, 
 $(X,G)$ is not stably conjugate 
to $(\PP^2,G)$ for any action of $G$ on $\PP^2$. 
\end{say}

%
 
 \newcommand{\etalchar}[1]{$^{#1}$}
\def\cprime{$'$} \def\polhk#1{\setbox0=\hbox{#1}{\ooalign{\hidewidth
 \lower1.5ex\hbox{`}\hidewidth\crcr\unhbox0}}}

\end{document}